\documentclass[12pt]{amsart}
\usepackage{amsmath,amsthm,amssymb,xypic,longtable}
\usepackage[matrix,arrow,curve]{xy}
\theoremstyle{plain}  

\newtheorem{thm}{Theorem}[section]
\newtheorem{prop}[thm]{Proposition}
\newtheorem{lem}[thm]{Lemma}   
\newtheorem{cor}[thm]{Corollary}

\theoremstyle{definition}

\newtheorem{defn}[thm]{Definition}

\theoremstyle{remark}

\newtheorem{claim}[thm]{Claim}
\newtheorem{rem}[thm]{Remark}
\newtheorem{ex}[thm]{Example}

\newcommand{\bA}{\mathbb A}
\newcommand{\bC}{\mathbb C}
\newcommand{\bF}{\mathbb F}
\newcommand{\bN}{\mathbb N}
\newcommand{\bP}{\mathbb P}
\newcommand{\bQ}{\mathbb Q}
\newcommand{\bR}{\mathbb R}
\newcommand{\bZ}{\mathbb Z}
\newcommand{\bG}{\mathbb G}

\newcommand{\cF}{\mathcal F}
\newcommand{\cO}{\mathcal O}
\newcommand{\cT}{\mathcal T}
\newcommand{\cX}{\mathcal X}
\newcommand{\cY}{\mathcal Y}

\newcommand{\oX}{\overline{X}}
\newcommand{\oE}{\overline{E}}
\newcommand{\oF}{\overline{F}}
\newcommand{\oM}{\overline{M}}

\newcommand{\cC}{\mathcal C}
\newcommand{\cE}{\mathcal E}
\newcommand{\cN}{\mathcal N}

\DeclareMathOperator{\Hom}{Hom}
\DeclareMathOperator{\Pic}{Pic}

\DeclareMathOperator{\Cl}{Cl}
\DeclareMathOperator{\Spec}{Spec}
\DeclareMathOperator{\Sing}{Sing}
\DeclareMathOperator{\Diff}{Diff}
\DeclareMathOperator{\Supp}{Supp}
\DeclareMathOperator{\Ext}{Ext}
\DeclareMathOperator{\End}{End}
\newcommand{\cHom}{\mathcal{H}om}
\newcommand{\cExt}{\mathcal{E}xt}

\DeclareMathOperator{\bmu}{\boldsymbol{\mu}}

\author{Paul Hacking}
\thanks{ The first author was partially supported by NSF grant DMS-0650052.
The second author was partially supported by grants RFBR no. 
\ 08-01-00395-a, 06-01-72017-MHTI-a and CRDF-RUM, no. 1-2692-MO-05.\\
\indent MSC classes: 14J10, 14E30}
\address{
Paul Hacking,  Department of Mathematics, University of Washington, Box 354350, Seattle, WA~98195}
\email{hacking@math.washington.edu} 
\author{
Yuri Prokhorov}

\address{
Yuri Prokhorov, Department of Higher Algebra, Faculty of Mathematics and Mechanics,  Moscow State Lomonosov University, Vorobievy Gory, Moscow, 119 899, RUSSIA}
\email{prokhoro@mech.math.msu.su}

\title{Smoothable del Pezzo surfaces with quotient singularities}
\date{}
\begin{document}
\maketitle

\section{Introduction}

We give a complete classification of del Pezzo surfaces with quotient singularities and Picard rank $1$ which admit a 
$\bQ$-Gorenstein smoothing. This solves \cite[Problem~28]{K2} in the case that the canonical class is negative.

Let $X$ be a normal surface with quotient singularities.
We say $X$ admits a \emph{$\bQ$-Gorenstein smoothing} if there exists a deformation $\cX/(0 \in T)$ of $X$ over a smooth curve germ such that the general fibre is smooth and $K_{\cX}$ is $\bQ$-Cartier.
(The requirement that $K_{\cX}$ be $\bQ$-Cartier is natural from the point of view of the minimal model program and is important in moduli problems, cf. \cite[5.4]{KSB}. It is automatically satisfied if $X$ is Gorenstein.)
We say $X'$ is a $\bQ$-Gorenstein deformation of $X$ if there exists a deformation $\cX/(0 \in T$) of $X$ over a smooth curve germ 
such that $\cX_t$ is isomorphic to $X'$ for all $t \neq 0$ and $K_{\cX}$ is $\bQ$-Cartier.

\begin{thm}\label{thmintro}
Let $X$ be a projective surface with quotient singularities such that $-K_X$ is ample and $\rho(X)=1$.
If $X$ admits a $\bQ$-Gorenstein smoothing then $X$ is either a $\bQ$-Gorenstein deformation of a toric surface with the same properties or one of the sporadic surfaces listed in Ex.~\ref{sporadic}.
\end{thm}

There are $14$ infinite families of toric examples, see Thm.~\ref{toric}. 
The surfaces in each family correspond to solutions of a Markov-type equation.
The solutions of the (original) Markov equation 
$$a^2+b^2+c^2=3abc$$
correspond to the vertices of an infinite tree such that each vertex has degree $3$.
Here two vertices are joined by an edge if they are related by a so called mutation of the form 
$$(a,b,c) \mapsto (a,b,3ab-c).$$
The solutions of the other equations  are described similarly.

Given one of the toric surfaces $Y$, the $\bQ$-Gorenstein deformations of $Y$ which preserve the Picard number are as follows.
First, there are no locally trivial deformations and no local-to-global obstructions to deformations.
Second, for each singularity $Q \in Y$, the deformation is either locally trivial or a deformation of a singularity of index $> 1$ 
to a Du Val singularity of type $A$, see Cor.~\ref{defspreservingrho}. 
Moreover, in the second case, the deformation is essentially unique (it is pulled back from a fixed one parameter deformation).

There are $20$ isolated sporadic surfaces and one family of sporadic surfaces parametrised by $\bA^1$, see Ex.~\ref{sporadic}.
Every sporadic surface has index $\le 2$. In particular, they occur in the list of Alexeev and Nikulin \cite{AN}.

In the case $K_X^2=9$ we obtain the following stronger result. This completely solves the problem studied by Manetti in \cite{M1}.
\begin{cor}\label{deg9}
Let $X$ be a projective surface with quotient singularities which admits a smoothing to the plane.
Then $X$ is a $\bQ$-Gorenstein deformation of a weighted projective plane $\bP(a^2,b^2,c^2)$,
where $(a,b,c)$ is a solution of the Markov equation.
\end{cor}

We note that a partial classification of the surfaces with $K_X^2 \ge 5$ was obtained by Manetti \cite{M1},\cite{M2}.

As a consequence of our techniques we verify a particular case 
of Reid's general elephant conjecture (see, e.g., \cite{Alexeev-1994ge}).
\begin{thm}\label{thm-ge}
Let $f\colon V\to (0 \in T)$ be a del Pezzo fibration over the germ of a smooth curve.
That is, $V$ is a $3$-fold with terminal singularities, 
$f$ has connected fibres, $-K_V$ is ample over $T$, and $\rho(V/T)=1$.
Assume in addition that the special fibre is reduced and normal, 
and has only quotient singularities.
Then a general member $S\in |-K_V|$ is a normal surface with Du Val singularities.
\end{thm}

In the final section we connect our results with the theory of exceptional vector bundles on del Pezzo surfaces.
Roughly speaking, given a $\bQ$-Gorenstein smoothing of a del Pezzo surface $X$ with quotient singularities,
there are exceptional vector bundles on the smooth fibre which are analogous to vanishing cycles.

\medskip
\noindent
\textbf{Notation}:
Throughout this paper, we work over the field $k=\bC$ of complex numbers.
The symbol $\bmu_n$ denotes the group of $n$th roots of unity.

\medskip

\noindent
\textbf{Acknowledgements}:
We thank I.~Dolgachev and J.~Koll\'ar for useful discussions.
We also thank J.~Stevens for pointing out the reference \cite{BC}.

\section{$T$-singularities} \label{secT}

$T$-singularities are by definition the quotient singularities of dimension $2$ which admit a $\bQ$-Gorenstein smoothing.
We recall the classification of $T$-singularities from \cite[Sec.~3]{KSB} and establish some basic results.

\subsection{$\bQ$-Gorenstein deformations}\label{QG}

We recall the definition and basic properties of $\bQ$-Gorenstein deformations of surfaces over an arbitrary base $S$ 
from \cite[Sec.~3]{H}. This definition was originally proposed by Koll\'ar \cite{K1}.

We first recall the notion of the canonical covering (or index one cover) of a $\bQ$-Gorenstein singularity.
Let $P \in X$ be a normal singularity such that $K_X$ is $\bQ$-Cartier. Let $n \in \bN$ be
the least integer such that $nK_X$ is Cartier, the \emph{index} of $P \in X$.
The \emph{canonical covering} $p : Z \rightarrow X$ of $P \in X$ is a Galois cover of $X$ with group $\bmu_n$, 
such that $Z$ is Gorenstein and $p$ is \'{e}tale in codimension $1$. Explicitly, 
$$Z= \underline{\Spec}_X(\cO_X \oplus \cO_X(K_X) \oplus \cdots \oplus \cO_X((n-1)K_X))$$ 
where the multiplication in $\cO_Z$ is given by fixing an isomorphism $\cO_X(nK_X) \stackrel{\sim}{\rightarrow} \cO_X$.
The canonical covering is uniquely determined up to isomorphism (assuming we work \'etale locally at $P \in X$).

\begin{defn}
Let $X$ be a normal surface such that $K_X$ is $\bQ$-Cartier.
We say that a deformation $\cX/(0 \in S)$ of $X$ is \emph{$\bQ$-Gorenstein} if, at every point $P \in X$,
$\cX/S$ is induced by an equivariant deformation of the canonical covering of $P \in X$. 
\end{defn}

\begin{rem}
Let $\omega_{\cX/S}$ denote the relative dualising sheaf of $\cX/S$.
Let $i \colon \cX^0 \subset \cX$ denote the Gorenstein locus of $\cX/S$ (i.e., the open locus where
$\omega_{\cX/S}$ is invertible).
For $M \in \bZ$, define $\omega_{\cX/S}^{[M]} = i_*\omega_{\cX^0/S}^{\otimes M}$.
Then $\cX/S$ is $\bQ$-Gorenstein iff $\omega_{\cX/S}^{[M]}$ commutes with base change for all $M \in \bZ$,
that is, for all $ f \colon T \rightarrow S$, the natural maps 
$$f^*\omega_{\cX/S}^{[M]} \rightarrow \omega_{\cX \times_S T /T}^{[M]}$$
are isomorphisms.
Moreover, in this case, $\omega_{\cX/S}^{[M]}$ is flat over $S$ for all $M \in \bZ$.
\end{rem}

\begin{rem}
If $\cX/S$ is $\bQ$-Gorenstein then $\omega_{\cX/S}^{[N]}$ is invertible for some $N \in \bN$.
(More precisely, for each $P \in X$, $\omega_{\cX/S}^{[N]}$ is invertible at $P \in \cX$ iff $\omega_X^{[N]}$ is invertible at 
$P \in X$ \cite[Lem.~3.3]{H}.) Conversely, if $S$ is a smooth curve, every fibre of $\cX/S$ has only quotient singularities,
and $\omega_{\cX/S}^{[N]}$ is invertible for some $N \in \bN$, then $\cX/S$ is $\bQ$-Gorenstein (cf. \cite[Lem.~3.4]{H}).
\end{rem}

The data of canonical coverings everywhere locally on $X$ defines a Deligne--Mumford stack $\mathfrak{X}$ with coarse moduli space $X$,
the \emph{canonical covering stack} of $X$. 
Explicitly, let $P \in X$ be a point, $n$ the index of $P \in X$, and $V \rightarrow U$ a canonical covering of a neighbourhood 
$U$ of $P \in X$. Then $\mathfrak{X}|_U$  is isomorphic to $[V/\bmu_n]$ over $U=V / \bmu_n$.
 
The deformations of the stack $\mathfrak{X}$ are naturally identified with the $\bQ$-Gorenstein deformations of $X$ (by passing to the coarse moduli space), see \cite[Prop.~3.7]{H}.
This is a useful point of view for studying global questions.

\begin{prop} \label{qtunobsQG}
Quotient singularities of dimension $2$ have unobstructed $\bQ$-Gorenstein deformations.
\end{prop}
\begin{proof}
The canonical covering of a quotient singularity is a Du~Val singularity, and so in particular a hypersurface singularity.
The $\bQ$-Gorenstein deformations are given by the equivariant deformations of the canonical covering,
so they are unobstructed. 
\end{proof}

\begin{rem} \label{obsTsing}
Quotient singularities typically have obstructed deformations.
For example, the deformation space of $\frac{1}{4}(1,1)$ has two smooth irreducible components of dimensions $1$ and 
$3$ which meet transversely at the origin \cite[8.2]{Pi}. The $1$-dimensional component is the $\bQ$-Gorenstein deformation space. 
See \cite[Thm.~3.9]{KSB} for a description of the components of the deformation space of an arbitrary quotient singularity.
\end{rem}

\subsection{Definition and classification of $T$-singularities}\label{Tdefclass}

\begin{defn}[{\cite[Def.~3.7]{KSB}}] \label{defT}
Let $P \in X$ be a quotient singularity of dimension $2$.
We say $P \in X $ is a \emph{$T$-singularity} if it admits a $\bQ$-Gorenstein smoothing.
That is, there exists a $\bQ$-Gorenstein deformation of $P \in X$ over a smooth curve germ
such that the general fibre is smooth.
\end{defn}

For  $n,a \in \bN$ with $(a,n)=1$, let $\frac{1}{n}(1,a)$ denote the cyclic quotient singularity 
$(0 \in \bA^2_{u,v} / \bmu_n)$ given by 
$$\bmu_n \ni \zeta \colon (u,v) \mapsto (\zeta u, \zeta^a v).$$

The following result is due to J.~Wahl \cite[5.9.1]{W2}, \cite[Props.~5.7,5.9]{LW}. 
It was proved by a different method in \cite[Prop.~3.10]{KSB}.

\begin{prop}\label{Tclass}
A $T$-singularity is either a Du~Val singularity or a cyclic quotient singularity of the form 
$\frac{1}{dn^2}(1,dna-1)$ for some $d,n,a \in \bN$ with $(a,n)=1$.
\end{prop}

The singularity $\frac{1}{dn^2}(1,dna-1)$ has index $n$ and canonical covering $\frac{1}{dn}(1,-1)$, the Du~Val singularity of type
$A_{dn-1}$.
We have an identification
$$\frac{1}{dn}(1,-1)=(xy=z^{dn}) \subset \bA^3_{x,y,z},$$ 
where $x=u^{dn}$, $y=v^{dn}$, and $z=uv$.
Taking the quotient by $\bmu_n$ we obtain
$$\frac{1}{dn^2}(1,dna-1) = (xy=z^{dn}) \subset \frac{1}{n}(1,-1,a).$$ 
Hence a $\bQ$-Gorenstein smoothing is given by
$$(xy=z^{dn}+t) \subset \frac{1}{n}(1,-1,a) \times \bA^1_t.$$
More generally, a versal $\bQ$-Gorenstein deformation of $\frac{1}{dn^2}(1,dna-1)$ is given by 
\begin{equation} \label{versalQGdef}
(xy=z^{dn}+t_{d-1}z^{(d-1)n}+\cdots+t_0) \subset \frac{1}{n}(1,-1,a) \times \bA^1_{t_0,\ldots,t_{d-1}}.
\end{equation}
We call a $T$-singularity of the form $\frac{1}{dn^2}(1,dna-1)$ a $T_d$-\emph{singularity}.

\begin{prop} \label{Tdeformation}
Let $(P \in \cX)/(0 \in S)$ be a $\bQ$-Gorenstein deformation of $\frac{1}{dn^2}(1,dna-1)$.
Then the possible singularities of a fibre of $\cX/S$ are as follows: 
either $A_{e_1-1},\ldots,A_{e_s-1}$ or \mbox{$\frac{1}{e_1n^2}(1,e_1na-1)$}, $A_{e_2-1},\ldots, A_{e_s-1}$,
where $e_1,\ldots,e_s$ is a partition of $d$.
\end{prop}
\begin{proof}
The family $\cX/S$ is pulled back from the versal $\bQ$-Gorenstein deformation (\ref{versalQGdef}).
Hence each fibre of $\cX/S$ has the form 
$$(xy=z^{dn}+a_{d-1}z^{(d-1)n}+\cdots+a_0) \subset \frac{1}{n}(1,-1,a)$$
for some $a_0,\ldots,a_{d-1} \in k$.
Write $$z^{dn}+a_{d-1}z^{(d-1)n}+\cdots+a_0= \prod (z^n-\gamma_i)^{e_i}$$ where the $\gamma_i$ are distinct.
Then the fibre has singularities as described in the statement (the second case occurs if $\gamma_i = 0$ for some $i$).
\end{proof}

\subsection{Noether's formula}

For $P \in X$ a $T$-singularity, let $M$ be the Milnor fibre of a $\bQ$-Gorenstein smoothing.
Thus $(M,\partial M)$ is a smooth $4$-manifold with boundary, and is uniquely determined by $P \in X$
since the $\bQ$-Gorenstein deformation space of $P \in X$ is smooth. Let $\mu_P= b_2(M)$, the \emph{Milnor number}.

\begin{lem}\cite[Sec.~3]{M1}\label{milnornumber}
If $P \in X$ is a Du Val singularity of type $A_r$, $D_r$, or $E_r$, then $\mu_P=r$.
If $P \in X$ is of type $\frac{1}{dn^2}(1,dna-1)$ then $\mu_P=d-1$. 
\end{lem}

\begin{rem}
If $M$ is the Milnor fibre of a smoothing of a normal surface singularity $P \in X$ then $M$ has the homotopy type of a 
CW complex of real dimension $2$ by Morse theory and $b_1(M)=0$ \cite{GS}. In particular $e(M)=1+\mu_P$.
\end{rem}

\begin{prop} \label{noether}
Let $X$ be a projective surface with $T$-singularities. Then 
$$K_X^2 + e(X) + \sum_{P \in \Sing X} \mu_P = 12\chi(\cO_X)$$
where $e(X)$ denotes the topological Euler characteristic 
and $\mu_P$ is the Milnor number. 

In particular, if $X$ is rational, then
$$K_X^2+\rho(X) + \sum_{P \in \Sing X} \mu_P = 10.$$
\end{prop}
\begin{proof}
For $X$ a normal surface with quotient singularities there is a singular Noether formula
$$K_X^2+e(X)+\sum_P c_P=12\chi(\cO_X)$$
where the sum is over the singular points $P \in X$,
and the correction term $c_P$ depends only on the local analytic isomorphism type of the singularity $P \in X$.
(Indeed, let $\pi \colon \tilde{X} \rightarrow X$ be the minimal resolution of $X$ and $E_1,\ldots,E_n$ the exceptional curves.
Noether's formula on $\tilde{X}$ gives $K_{\tilde{X}}^2+e(\tilde{X})=12\chi(\cO_{\tilde{X}})$.
Write $K_{\tilde{X}}=\pi^*K_X+\sum a_iE_i = \pi^*K_X +A$.
Then $K_{\tilde{X}}^2=K_X^2+A^2$, $e(\tilde{X})=e(X)+n$ (by the Mayer--Vietoris sequence), 
and $\chi(\cO_{\tilde{X}})=\chi(\cO_X)$ (because $X$ has rational singularities). Hence $K_X^2+e(X)+(A^2+n)=12\chi(\cO_X)$.)

For each singularity $P \in X$, there exists a projective surface $Y$ with a unique singularity $Q \in Y$ which is isomorphic to 
$P \in X$, and a $\bQ$-Gorenstein smoothing $\cY/(0 \in T)$ of $Y$ over a smooth curve germ
(this is a special case of Looijenga's globalisation theorem \cite[App.]{L}).
Let $Y'$ be a general fibre of $\cY/T$, then $K_{Y'}^2+e(Y')=12\chi(\cO_{Y'})$.
Now $K_{Y'}^2=K_{Y}^2$, $\chi(\cO_{Y'})=\chi(\cO_Y)$, and $e(Y') = e(Y) + \mu_Q$ 
(because the Milnor fibre of the smoothing has Euler number $1+\mu_Q$).
Hence $K_Y^2+e(Y)+\mu_Q = 12\chi(\cO_Y)$. Thus the correction term to Noether's formula due to the $T$-singularity $Q \in Y$ is $\mu_Q$.
This gives the result.
\end{proof}

\begin{cor}\label{defspreservingrho}
Let $X$ be a projective surface with $T$-singularities and $X'$ a fibre of a $\bQ$-Gorenstein deformation $\cX/(0 \in T)$ of $X$ over a smooth curve germ.
Then $e(X)=e(X')$ iff at each singular point $P \in X$, the deformation is either locally trivial or a deformation of a 
$T_d$-singularity to an $A_{d-1}$ singularity.
\end{cor}  
\begin{proof}
This follows immediately from Props.~\ref{Tdeformation} and \ref{noether}.
\end{proof}

\subsection{Minimal resolutions of $T$-singularities}

Given a cyclic quotient singularity $\frac{1}{n}(1,a)$, let $[b_1,\ldots,b_r]$ be the expansion of $n/a$ as a 
Hirzebruch--Jung continued fraction \cite[p.~46]{F}.
Then the exceptional locus of the minimal resolution of $\frac{1}{n}(1,a)$ is a chain of smooth rational curves of 
self-intersection numbers $-b_1,\ldots,-b_r$. 
The strict transforms of the coordinate lines $(u=0)$ and $(v=0)$ intersect the right and left end components
of the chain respectively.
\begin{rem} \label{reverse}
Note that $[b_r,\ldots,b_1]$ corresponds to the same singularity as $[b_1,\ldots,b_r]$ with the roles of the coordinates 
$u$ and $v$ interchanged. Thus, if $[b_1,\ldots,b_r]=n/a$ then $[b_r,\ldots,b_1]=n/a'$ where $a'$ is the inverse of $a$ 
modulo $n$.
\end{rem}

We recall the description of the minimal resolution of the cyclic quotient singularities of class $T$.
Let a \emph{$T_d$-string} be a string $[b_1,\ldots,b_r]$ which corresponds to a $T_d$-singularity.
\begin{prop}\cite[Prop.~3.11]{KSB}, \cite[Thm.~17]{M1} \label{exclocusT_d}
\begin{enumerate}
\item $[4]$ is a $T_1$-string and, for $d \ge 2$, $[3,2,\ldots,2,3]$ (where there are $(d-2)$ $2$'s)
is a $T_d$-string.
\item If $[b_1,\ldots,b_r]$ is a $T_d$-string, then so are $[b_1+1,b_2,\ldots,b_r,2]$ and $[2,b_1,\ldots,$ $b_r+1]$.
\item For each $d$, all $T_d$-strings are obtained from the example in 
\textup{(1)} by iterating the steps in \textup{(2)}.
\end{enumerate}
\end{prop}

\section{Unobstructedness of deformations}

We prove that for a del Pezzo surface with $T$-singularities there are no local-to-global obstructions to deformations. 
Thus a del Pezzo surface with quotient singularities admits a $\bQ$-Gorenstein smoothing iff it has $T$-singularities.

\begin{lem} \label{RR}
Let $X$ be a projective surface such that $X$ has only $T$-singularities and $-K_X$ is nef and big. 
Then $$h^0(\cO_X(-nK_X))=1+\frac{1}{2}n(n+1)K_X^2$$ for $n \in \bZ_{\ge 0}$.
\end{lem}

\begin{proof}
For $X$ a projective normal surface with only quotient singularities and $D$ a Weil divisor on $X$,
we have a singular Riemann--Roch formula
$$\chi(\cO_X(D))=\chi(\cO_X)+ \frac{1}{2}D(D-K_X) + \sum c_P(D),$$
where the sum is over points $P \in X$ where the divisor $D$ is not Cartier 
and the correction term $c_P(D)$ depends only on the local analytic isomorphism type of the singularity $P \in X$ 
and the local analytic divisor class of $D$ at $P \in X$ \cite[1.2]{B}.
We prove that $c_P(mK_X)=0$ for $P \in X$ a $T$-singularity and $m \in \bZ$.
There exists a projective surface $Y$ with a unique singularity $Q \in Y$ isomorphic to $P \in X$
and a $\bQ$-Gorenstein smoothing $\cY/(0 \in T)$ of $Y$ over a smooth curve germ 
(by Looijenga's globalisation theorem \cite[App.]{L}).
We compute that $c_Q(mK_Y)=0$ for all $m \in \bZ$ by comparing the Riemann--Roch formulae on $Y$ and a general fibre $Y'$ of $\cY/T$. 
The Riemann--Roch formula for the line bundle $\cO_{Y'}(mK_{Y'})$ on $Y'$ gives 
$\chi(\cO_{Y'}(mK_{Y'}))=\chi(\cO_{Y'})+ \frac{1}{2}m(m-1)K^2_{Y'}$.
Now $\chi(\cO_{Y'})=\chi(\cO_Y)$, $\chi(\cO_{Y'}(mK_{Y'}))=\chi(\cO_Y(mK_Y))$ 
(note that $\omega_{\cY/T}^{[m]}$ is flat over $T$ and commutes with base change because $\cY/T$ is $\bQ$-Gorenstein), 
and $K_{Y'}^2=K_Y^2$.  
Hence $\chi(\cO_Y(mK_{Y}))=\chi(\cO_{Y})+\frac{1}{2}m(m-1)K_Y^2$, i.e., $c_Q(mK_Y)=0$.

Now suppose that $X$ has only $T$-singularities and $-K_X$ is nef and big as in the statement.
Then $\chi(\cO_X(mK_X))=\chi(\cO_X)+\frac{1}{2}m(m-1)K_X^2$ for $m \in \bZ$
and $H^i(\cO_X(-nK_X))=H^i(\cO_X)=0$ for $i>0$ and $n \ge 0$ by Kawamata--Viehweg vanishing. 
Hence $h^0(\cO_X(-nK_X))=1+\frac{1}{2}n(n+1)K_X^2$, as required.
\end{proof}

\begin{prop} \label{unobs}
Let $X$ be a projective surface such that $X$ has only $T$-singularities and 
$-K_X$ is nef and big. 
Then there are no local-to-global obstructions to deformations of $X$.
In particular, $X$ admits a $\bQ$-Gorenstein smoothing.
Moreover $X$ has unobstructed $\bQ$-Gorenstein deformations.
\end{prop}

\begin{proof} (cf. \cite[Pf. of Thm.21]{M1})
The local-to-global obstructions to deformations of $X$ 
lie in $H^2(T_X)$, where $T_X =\cHom(\Omega_X,\cO_X)$ is the tangent sheaf of $X$.
This follows from either a direct cocycle computation (cf. \cite[Prop.~6.4]{W2}) or the theory of the cotangent complex
\cite[2.1.2.3]{I}.
We have $H^2(T_X)=\Hom(T_X,\cO_X(K_X))^*$ by Serre duality.
Since $H^0(-K_X) \neq 0$ by Lem.~\ref{RR}, we have an inclusion
$$\Hom(T_X,\cO_X(K_X)) \subset \Hom(T_X,\cO_X) = H^0(\Omega_X^{\vee\vee}).$$
Here $\Omega_X^{\vee\vee}$ is the double dual or reflexive hull of $\Omega_X$.
By Steenbrink's Hodge theory for orbifolds \cite[Sec.~1]{S}, we have $h^0(\Omega_X^{\vee\vee})=h^1(\cO_X)=0$.
Combining, we deduce $H^2(T_X)=0$. So there are no local-to-global obstructions for deformations of $X$.
$T$-singularities have unobstructed $\bQ$-Gorenstein deformations by Prop.~\ref{qtunobsQG}.
Hence $X$ has unobstructed $\bQ$-Gorenstein deformations.
\end{proof}

\begin{rem}
The surface $X$ has obstructed deformations in general because $T$-singularities have obstructed deformations, see Rem.~\ref{obsTsing}.
\end{rem}

\begin{rem}
There may be local-to-global obstructions in positive characteristic.
For example, let $k$ be an algebraically closed field of characteristic $2$. 
Let $\tilde{X}$ be the blowup of $\bP^2_k$ in the set of $7$ points 
$\bP^2_{\bF_2} \subset \bP^2_k$.
Each line in $\bP^2_k$ which is defined over $\bF_2$ passes through $3$ of the points.
Let $\pi \colon \tilde{X} \rightarrow X$ be the contraction of the $7$ $(-2)$-curves on $\tilde{X}$ 
given by the strict transforms of these lines.
Then $X$ is a Gorenstein log del Pezzo surface with $7$ $A_1$ singularities.

We claim that there are local-to-global obstructions to deformations of $X$.
Suppose this is not the case.
Since $\tilde{X}$ is the blowup of $7$ distinct points in $\bP^2$, $4$ of which are in general position,
it is easy to see that $\tilde{X}$ has no infinitesimal automorphisms and its universal deformation space is smooth of 
dimension $6$ (the deformations of $\tilde{X}$ are given by moving the points we blow up).
Moreover, since $H^1(\cO_{\tilde{X}})=0$, there is a ``blowing down map'' from deformations of $\tilde{X}$ to deformations of $X$
\cite[Thm.~1.4(c)]{W1}.
In particular, at first order, we have a map $H^1(T_{\tilde{X}}) \rightarrow \Ext^1(\Omega_X,\cO_X)$. 

The $A_1$ singularity is the hypersurface singularity
$$(xy=z^2) \subset \bA^3_{x,y,z}$$
(this is true in any characteristic).
By our assumption, there exists a deformation $\cX/S$ of $X$ over $S=(0 \in \bA^7_{t_1,\ldots,t_7})$ inducing a
deformation of the $i$th singular point of the form
$$(xy=z(z+t_i)) \subset \bA^3_{x,y,z} \times S$$
(note that this deformation is non-trivial at first order in characteristic $2$).
There exists a simultaneous resolution of the family $\cX/S$ --- near each singular point of the special fibre, 
we blow up the locus $(z=x=0) \subset \cX$ to obtain a small resolution $f \colon \tilde{\cX} \rightarrow \cX$ 
which restricts to the minimal resolution of each fibre of $\cX/S$.
Hence the Kodaira--Spencer map 
$$T_S \otimes k(0) \rightarrow \Ext^1(\Omega_X,\cO_X)$$
for the family $\cX/S$ at $0 \in S$ factors through $H^1(T_{\tilde{X}})$.
But $h^1(T_{\tilde{X}})=6$ and the composition
$$T_S \otimes k(0) \rightarrow \Ext^1(\Omega_X,\cO_X) \rightarrow H^0(\cExt^1(\Omega_X,\cO_X))$$
(the Kodaira--Spencer map for the induced deformation of the singularities) is injective by construction, a contradiction.
\end{rem}

\section{Toric surfaces}

\begin{thm} \label{toric}
The projective toric surfaces with $T$-singularities and Picard rank $1$ are as follows.
There are $14$ infinite families $(1),\ldots,(8.4)$ which we list in the tables below. 
In cases $(1),\ldots,(4)$, the surface $X$ is a weighted projective plane $\bP(w_0,w_1,w_2)$, and the weights $w_0,w_1,w_2$
are determined by a solution $(a,b,c)$ of a Markov-type equation.
In the remaining cases, the surface $X$ is a quotient of one of the above weighted projective planes $Y$ 
by $\bmu_e$ acting freely in codimension $1$.
The action is diagonal with weights $(m_0,m_1,m_2)$, i.e.,
$$\bmu_e \ni \zeta \colon (X_0,X_1,X_2) \mapsto (\zeta^{m_0} X_0, \zeta^{m_1} X_1, \zeta^{m_2} X_2)$$
where $X_0,X_1,X_2$ are homogeneous coordinates on $Y$.
We also record $K_X^2$ and the values of $d=\mu+1$ for the singularities of $X$.

\begin{eqnarray*} 
\renewcommand{\arraystretch}{1.5}
\begin{array}{|l|l|l|l|l|} 
\hline
X	& w_0,w_1,w_2 	& \mbox{\rm Markov-type equation}  & K_X^2		& d \\
\hline
(1)	& a^2,b^2,c^2	& a^2+b^2+c^2=3abc 	& 9		& 1, 1, 1 \\
(2)	& a^2,b^2,2c^2	& a^2+b^2+2c^2=4abc 	& 8		& 1, 1, 2 \\        	       
(3)	& a^2,2b^2,3c^2	& a^2+2b^2+3c^2=6abc	& 6		& 1, 2, 3 \\
(4)	& a^2,b^2,5c^2	& a^2+b^2+5c^2=5abc	& 5		& 1, 1, 5 \\
\hline
\end{array}
\end{eqnarray*}
\begin{eqnarray*} 
\renewcommand{\arraystretch}{1.5}
\begin{array}{|l|l|l|l|l|l|} 
\hline
X	& Y \quad & e \quad		& m_0,m_1,m_2      	& K_X^2		& d \\
\hline
(5)	& (2)	& 2	 	& 0, 1, -1		& 4	 	& 2, 2, 4 \\
(6.1)   & (1)	& 3		& 0, 1, -1		& 3		& 3, 3, 3 \\
(6.2)	& (3)	& 2		& 0, 1, -1		& 3		& 1, 2, 6 \\
(7.1)	& (2)	& 4		& 0, 1, 1		& 2		& 1, 1, 8 \\
(7.2)	& (2)	& 4		& 0, 1, -1		& 2 		& 2, 4, 4 \\
(7.3)	& (3)	& 3		& 0, 1, -1		& 2		& 1, 3, 6 \\
(8.1)	& (1)	& 9		& 0, 1, -1		& 1		& 1, 1, 9 \\
(8.2)	& (2)	& 8		& 0, 1, -1		& 1		& 1, 2, 8 \\
(8.3)	& (3)	& 6		& 0, 1, -1		& 1		& 2, 3, 6 \\
(8.4)	& (4)	& 5		& 0, 1, -1		& 1		& 1, 5, 5 \\
\hline
\end{array}
\end{eqnarray*}
\end{thm}

\begin{rem}
With notation as above, 
let $X^0 \subset X$ be the smooth locus and $p^0 \colon Y^0 \rightarrow X^0$ the restriction of the cover $Y \rightarrow X$.
Then $p^0$ is the universal cover of  $X^0$. In particular $\pi_1(X^0)$ is cyclic of order $e$.
\end{rem}

The solutions of the Markov-type equations in Thm.~\ref{toric} may be described as follows \cite[3.7]{KN}. 
We say a solution $(a,b,c)$ is \emph{minimal} if $a+b+c$ is minimal.
The equations (1),(2),(3) have a unique minimal solution $(1,1,1)$, and (4) has minimal solutions
$(1,2,1)$ and $(2,1,1)$. 
Given one solution, we obtain another by regarding the equation as a quadratic in one of the variables, $c$ (say), 
and replacing $c$ by the other root. 
Explicitly, if the equation is $\alpha a^2 +\beta b^2 + \gamma c^2 =\lambda abc$, then
\begin{equation} \label{mutation}
(a,b,c) \mapsto (a,b, \frac{\lambda}{\gamma}ab-c).
\end{equation}
This process is called a \emph{mutation}.
Every solution is obtained from a minimal solution by a sequence of mutations.

For each equation, we define an infinite graph $\Gamma$ such that the vertices are labelled by the solutions and two vertices 
are joined by an edge if they are related by a mutation. For equation (1), $\Gamma$ is an infinite tree such that
each vertex has degree $3$, and there is an action of $S_3$ on $\Gamma$ given by permuting the variables $a,b,c$.
The other cases are similar, see \cite[3.8]{KN} for details.

\begin{proof}[Proof of Theorem~\ref{toric}]
Let $X$ be a projective toric surface such that $X$ has only $T$-singularites and $\rho(X)=1$.
The surface $X$ is given by a complete fan $\Sigma$ in $N_{\bR} \simeq \bR^2$, 
where $N \simeq \bZ^2$ is the group of $1$-parameter subgroups of the torus. 
The fan $\Sigma$ has $3$ rays because $\rho(X)=1$.
Let $v_0,v_1,v_2 \in N$ be the minimal generators of the rays.
There is a unique relation $$w_0v_0+w_1v_1+w_2v_2=0$$ where $w_0,w_1,w_2 \in \bN$ are pairwise coprime.
Let $N_Y \subseteq N$ denote the subgroup generated by $v_0,v_1,v_2$.
Let $p \colon Y \rightarrow X$ be the finite toric morphism corresponding to the inclusion $N_Y \subseteq N$. 
Then $Y$ is isomorphic to the weighted projective plane $\bP(w_0,w_1,w_2)$
and $p$ is a cyclic cover of degree $e=|N/N_Y|$ which is \'{e}tale over the smooth locus $X^0 \subset X$.
The surface $Y$ has only $T$-singularities because a cover of a $T$-singularity which is \'{e}tale in codimension $1$
is again a $T$-singularity (this follows easily from the classification of $T$-singularities).

The surface $X$ has $3$ cyclic quotient singularities of class $T$.
Let the singularities of $X$ be $\frac{1}{d_in_i^2}(1,d_in_ia_i-1)$ for $i=0,1,2$.
Then 
\begin{equation}
d_0+d_1+d_2+K_X^2=12
\end{equation} 
by Prop.~\ref{noether}.
The singularities of $X$ are quotients of the singularities 
$\frac{1}{w_0}(w_1,w_2)$, $\frac{1}{w_1}(w_0,w_2)$, $\frac{1}{w_2}(w_0,w_1)$ of $Y$
by $\bmu_e$. Hence $d_in_i^2=ew_i$. Also $K_Y^2=eK_X^2$ because $p \colon Y \rightarrow X$ has degree $e$ and is \'etale in codimension $1$.
Let $H$ be the ample generator of the class group of $Y$.
Then $K_Y \sim -(w_0+w_1+w_2)H$, and $H^2=\frac{1}{w_0w_1w_2}$.
We deduce that
\begin{equation} \label{MEpf}
d_0n_0^2+d_1n_1^2+d_2n_2^2=\sqrt{K_X^2d_0d_1d_2}\cdot n_0n_1n_2.
\end{equation}
In particular 
$$\sqrt{K_X^2d_0d_1d_2} = \sqrt{(12-\sum d_i)d_0d_1d_2} \in \bZ$$
We compute all triples $d=(d_0,d_1,d_2)$ satisfying this condition.
They are as listed in the last column of the tables above.

We first treat the cases $d=(1,1,1)$, $(1,1,2)$, $(1,2,3)$, and $(1,1,5)$. 
These are the cases for which $K_X^2 \ge 5$. Since $K_Y^2=e K_X^2 \le 9$ by Prop.~\ref{noether} we deduce that $e=1$.  
Thus $X$ is isomorphic to a weighted projective plane.
The weights $d_in_i^2$ are determined by the solution $(n_0,n_1,n_2)$ of (\ref{MEpf}), which is the Markov-type equation given in the statement.
Conversely, we check that for any solution of (\ref{MEpf}) the weighted projective plane $X=\bP(d_0n_0^2,d_1n_1^2,d_2n_2^2)$
has $T$ singularities and the expected value of $d$. 
We use the description of the solutions of (\ref{MEpf}) given above.
We write $\lambda=\sqrt{K_X^2d_0d_1d_2}$, and note that $d_0d_1d_2$ divides $\lambda$ in each case. 
By induction using (\ref{mutation}) we find that $n_0,n_1,n_2$ are pairwise coprime and $\gcd(n_i,\frac{\lambda}{d_i})=1$ for each $i$.
In particular, the $d_in_i^2$ are pairwise coprime. 
Now consider the singularity $\frac{1}{d_0n_0^2}(d_1n_1^2,d_2n_2^2)$.
We have $$d_1n_1^2+d_2n_2^2 = \lambda n_0n_1n_2 \mod d_0n_0^2$$ by (\ref{MEpf}), and 
so $\gcd(d_1n_1^2+d_2n_2^2,d_0n_0^2)=d_0n_0$ because $\gcd(\frac{\lambda}{d_0}n_1n_2,n_0)=1$.
Thus this singularity is of type $T_{d_0}$.

For the remaining values of $d$, we determine the degree $e$ of the cover $p \colon Y \rightarrow X$ as follows.
We have $e=\gcd(d_0n_0^2,d_1n_1^2,d_2n_2^2)$.
By inspecting the equation $(\ref{MEpf})$ we find a factor of $e$, and, together with the inequality $eK_X^2=K_Y^2 \le 9$,
this is sufficient to determine $e$ in each case.
For example, let $d=(1,2,8)$. Then we find that $n_0$ is divisible by $4$ and $n_1$ is even, so $e$ is divisible by $8$, hence
equal to $8$. In each case we have $K_Y^2 \ge 5$, so $Y$ is one of the surfaces classified above.

We now classify the possible actions of $\bmu_e$ on the covering surface $Y$.
We have $Y=\bP(d_0n_0^2,d_1n_1^2,d_2n_2^2)$ where $d=d_Y=(1,1,1),(1,1,2),(1,2,3),$ or $(1,1,5)$, and $(n_0,n_1,n_2)$ is a solution of
(\ref{MEpf}).
The action is given by
$$\bmu_e \ni \zeta \colon (X_0,X_1,X_2) \mapsto (\zeta^{m_0} X_0,\zeta^{m_1} X_1,\zeta^{m_2} X_2)$$
where $X_0,X_1,X_2$ are the homogeneous coordinates on the weighted projective plane $Y$.
In each case $d_0n_0^2=n_0^2$ is coprime to $e$. So we may assume that $m_0=0$.
We may also assume that $m_1=1$ (because the action is free in codimension $1$).
Consider the singularity $P_0 \in X$ below $(1:0:0)\in Y$. This singularity admits a covering by $\frac{1}{e}(1,m_2)$
(which is \'etale in codimension $1$). Hence $\frac{1}{e}(1,m_2)$ is a $T$-singularity.
If $e$ is square-free, it follows that $m_2=-1$.
If $e=4$, then $m_2=\pm 1$.
If $e=8$ then $d_Y=(1,1,2)$ and $d_X=(1,2,8)$, so we may assume that $P_0 \in X$ is a $T_8$-singularity 
(note that a $\bmu_8$-quotient of a $T_2$-singularity cannot be a $T_8$-singularity). 
Thus $P_0 \in X$ is covered by $\frac{1}{8}(1,-1)$ and so $m_2=-1$.
Similarly if $e=9$ then $d_Y=(1,1,1)$ and $d_X=(1,1,9)$, so we may assume that $P_0 \in X$ is a $T_9$-singularity, and $m_2=-1$.
This gives the list of group actions above.
Finally, it remains to check that for each such quotient $X=Y/\bmu_e$, the surface $X$ has $T$-singularities with the expected values
of $d$. This is a straightforward toric calculation, so we omit it.
\end{proof}

\section{Surfaces with a $D$ or $E$ singularity}

A \emph{log del Pezzo surface} is a normal projective surface $X$ such that $X$ has only quotient singularities and $-K_X$ is ample.

\begin{thm}\label{DorE}
Let $X$ be a log del Pezzo surface such that $\rho(X)=1$, and assume that $\dim |-K_X| \ge 1$.
\begin{enumerate}
\item If $X$ has a Du Val singularity 
of type $E$ then $K_X$ is Cartier.
\item If $X$ has a Du Val singularity of type $D$
then either $K_X$ is Cartier or there is a unique non Du~Val singularity 
of type $\frac{1}{m}(1,1)$ for some $m \ge 3$.
\end{enumerate}
Moreover, in both cases, a general member of $|-K_X|$ is irreducible and does not pass through the Du~Val singularities.
\end{thm}
\begin{proof}
Assume that $X$ has a $D$ or $E$ singularity $P \in X$ and $K_X$ is \emph{not} Cartier.
Let $\nu \colon \hat{X} \rightarrow X$ be the minimal resolution of the non Du~Val singularities of $X$ and write $\hat{P}=\nu^{-1}(P)$.
So $\hat{X}$ has only Du~Val singularities and $\hat{P} \in \hat{X}$ is a $D$ or $E$ singularity.
Let $\{ E_i \}$ be the exceptional curves of $\nu$ and write $E = \sum E_i$.

Write $|-K_{\hat{X}}|=|M|+F$ where $F$ is the fixed part and $M$ is general in $|M|$.
We have an equality
$$K_{\hat{X}}=\nu^*K_X+ \sum a_iE_i$$
where $a_i < 0$ for all $i$ because $\nu$ is minimal and we only resolve the non Du~Val singularities \cite[Lem.~3.41]{KM}.
Hence $\dim |-K_{\hat{X}}| = \dim |-K_X|$ and $F \ge E$.
 
We run the minimal model program on $\hat{X}$.
We obtain a birational morphism $\phi \colon \hat{X} \rightarrow \oX$ such that $\oX$ has Du~Val singularities and exactly one of the following holds.
\begin{enumerate}
\item $K_{\oX}$ is nef.
\item $\rho(\oX)=2$ and there is a fibration $\psi \colon \oX \rightarrow \bP^1$ with $K_{\oX} \cdot f < 0$ for $f$ a fibre.
\item $\rho(\oX)=1$ and $-K_{\oX}$ is ample.
\end{enumerate}
Clearly $K_{\oX}$ is not nef because $\dim|-K_{\oX}| \ge \dim |-K_{\hat{X}}| \ge 1$.

In the minimal model program for surfaces with Du~Val singularities, the birational extremal contractions are weighted blowups $f \colon X \rightarrow Y$
of a smooth point $P \in Y$ with weights $(1,n)$ for some $n \in \bN$.
In particular the exceptional divisor $E \subset X$ is a smooth rational curve and passes through a unique singularity of $X$ which is of type 
$\frac{1}{n}(1,-1)=A_{n-1}$.
See \cite[Lem.~3.3]{KMcK}.

Therefore, the birational morphism $\phi$ is an isomorphism near the $D$ or $E$ singularity $\hat{P} \in \hat{X}$
and $\oE := \phi_* E$ is contained in the smooth locus of $\oX$.
Note also that $\oE \neq 0$ because $\rho(X)=1$ and $X$ has a non Du~Val singularity.

Suppose first we are in case (3). We have $-K_{\oX} \sim \oM+\oF$ where $\oM:=\phi_*M$ is mobile and $\oF:=\phi_*F \ge \oE$.
In particular, $\Pic(\oX)$ is not generated by $-K_{\oX}$ because $\oM + \oF > \oE$ and $\oE$ is Cartier.
Hence $\oX$ is isomorphic to $\bP^2$ or $\bP(1,1,2)$ by the classification of Gorenstein log del Pezzo surfaces 
\cite{D}. (Indeed, if $Y$ is a Gorenstein del Pezzo surface, let $f \colon \tilde{Y} \rightarrow Y$ be the minimal resolution.
Then either $Y$ is isomorphic to $\bP^2$ or $\bP(1,1,2)$, or $\tilde{Y}$ is obtained from $\bP^2$ by a sequence of blowups.
In the last case, let $C \subset \tilde{Y}$ be a $(-1)$-curve.
Then 
$$K_Y \cdot f_*C = f^*K_{Y} \cdot C = K_{\tilde{Y}} \cdot C = -1 .$$
It follows that $-K_Y$ is a generator of $\Pic Y$ if $\rho(Y)=1$.)
So $\oX$ does not have a $D$ or $E$ singularity, a contradiction.

So we are in case (2). Write $p = \psi \circ \phi \colon \hat{X} \rightarrow \bP^1$.
The divisor $E$ has a $p$-horizontal component, say $E_1$ (because $\rho(X)=1$ so there does not exist a morphism $X \rightarrow \bP^1$).
If $f$ is a general fibre of $p$ then 
$$2=-K_{\hat{X}} \cdot f \ge E_1 \cdot f \ge 1.$$
If $E_1 \cdot f =1$ then all fibres of $\psi$ are reduced (because $\oE_1$ is contained in the smooth locus of $\oX$), 
so $\oX$ is smooth \cite[Lem.~11.5.2]{KMcK}, a contradiction.
So $E_1 \cdot f =2$.
Then $(M+(F-E_1))\cdot f =0$, so $M$ and $F-E_1$ are $p$-vertical.
In particular $M$ is basepoint free and $E_1$ has coefficient $1$ in $F$.
Since
$$
2 \ge 2-2p_a(E_1)=-(K_{\hat{X}}+E_1) \cdot E_1 = (M+(F-E_1)) \cdot E_1 \ge M \cdot E_1 \ge 2,
$$
we find $M \cdot E_1 = 2$ and $(F-E_1)\cdot E_1=0$.
Thus $M$ is a fibre of $\psi$ and the divisors $M+E_1$ and $F-E_1$ have disjoint support.
But $M+F \sim -K_{\hat{X}}$ is connected because 
$$
H^1(\cO_{\hat{X}}(-M-F))=H^1(K_{\hat{X}})=H^1(\cO_{\hat{X}})^*=0.
$$
Hence $F=E=E_1$. In particular, $X$ has a unique non Du~Val singularity of type $\frac{1}{m}(1,1)$ (where $E_1^2=-m$).
Also, a general member of $|-K_X|$ is irreducible and does not pass through any Du~Val singularities.
Finally $\oX$ does not have a singularity of type $E$ by the classification of fibres of $\bP^1$ fibrations with Du~Val singularities 
\cite[Lem.~11.5.12]{KMcK}. So $X$ does not have an $E$ singularity.

If $K_X$ is Cartier then a general member of $|-K_X|$ is smooth and misses the singular points by \cite{D}.
\end{proof}

\section{Surfaces of index $\le 2$}
 
Alexeev and Nikulin classified log del Pezzo surfaces $X$ of index $\le 2$ \cite{AN}.
They prove that $X$ is a $\bZ/2\bZ$ quotient of a K3 surface and use the Torelli theorem for K3 surfaces
to obtain the classification.
In this section, we deduce the index $\le 2$ case of our main theorem from their result.

We note that the quotient singularities of index $\le 2$ are the Du Val singularities and the cyclic quotient singularities
of type $\frac{1}{4d}(1,2d-1)$, see \cite{AN}. In particular, they are $T$-singularities.

\begin{prop} \label{index2deformation}
Let $X$ be a log del Pezzo surface of index $\le 2$ such that $\rho(X)=1$.
Then exactly one of the following holds.
\begin{enumerate}
\item $X$ is a $\bQ$-Gorenstein deformation of a toric surface.
\item $X$ has either a $D$ singularity, an $E$ singularity,  or $\ge 4$ Du~Val singularities.
\end{enumerate}
\end{prop}
\begin{proof}
We first observe that the two conditions cannot both hold.
If $X$ is a $\bQ$-Gorenstein deformation of a toric surface $Y$, then necessarily $\rho(Y)=1$ and $Y$ has only $T$-singularities.
In particular, $Y$ has at most $3$ singularities.
Moreover, since the deformation preserves the Picard number, the only possible non-trivial deformation of a singularity of $Y$ is
a deformation of a $T_d$ singularity to a $A_{d-1}$ singularity by Cor.~\ref{defspreservingrho}.
Finally, note that $Y$ does not have a $D$ or $E$ singularity because $Y$ is toric.
Hence $X$ has at most $3$ singularities and does not have a $D$ or $E$ singularity.

We now use the classification of log del Pezzo surfaces of index $\le 2$ and Picard rank $1$ \cite[Thms.~4.2,4.3]{AN}.
We check that each such surface $X$ which does not satisfy condition (2) is a deformation of a toric surface $Y$.
By \cite{AN}, $X$ is determined up to isomorphism by its singularities.
So it suffices to exhibit a toric surface $Y$ such that $\rho(Y)=1$ and the singularities of $X$ are obtained from the singularities of 
$Y$ by a $\bQ$-Gorenstein deformation which preserves the Picard number. We list the surfaces $Y$ in the tables below.
\end{proof}

In the following tables, for each log del Pezzo surface $X$ of Picard rank $1$ and index $\le 2$ such that $X$ does not satisfy condition (2) of 
Prop.~\ref{index2deformation}, 
we exhibit a toric surface $Y$ such that $X$ is a $\bQ$-Gorenstein deformation of $Y$.
The tables treat the surfaces $X$ of index $1$ and $2$ respectively.
We give the number of the surface $X$ in the list of Alexeev and Nikulin \cite[p.~93--100]{AN}.
We use the description of the toric surfaces $Y$ given in Thm.~\ref{toric}.
We give the number of the infinite family to which $Y$ belongs and the solution $(a,b,c)$ of the Markov-type equation
corresponding to $Y$.
We record the value of $d= \mu + 1$ for each singularity in the last column of the table.

\par\medskip\noindent
\setlongtables\renewcommand{\arraystretch}{1.3}
\begin{longtable}{|l|l|l|l|l|}
\hline
$X$&$\Sing X$&$Y$&$\Sing Y$&$d$ \\
\hline
\endfirsthead
\hline
$X$&$\Sing X$&$Y$&$\Sing Y$&$d$ \\
\hline
\endhead
\hline
\endlastfoot
\hline
\endfoot
1 & &$(1),\, (1,1,1)$& & 1, 1, 1 \\
2 &$A_1$&$(2),\, (1,1,1)$&$A_1$& 1, 1, 2 \\
5 &$A_1, A_2$&$(3),\, (1,1,1)$&$A_1, A_2$& 1, 2, 3 \\
6 &$A_4$&$(4),\, (1,2,1)$&$\frac{1}{4}(1,1)$, $A_4$ & 1, 1, 5 \\
7b &$2A_1$, $A_3$&$(5),\, (1,1,1)$&$2A_1$, $A_3 $ & 2, 2, 4 \\
8b &$A_1$, $A_5$&$(6.2),\, (1,1,1)$&$\frac{1}{4}(1,1)$, $A_1$, $A_5$ & 1, 2, 6 \\
8c &$3A_2$&$(6.1),\, (1,1,1)$&$3A_2$ & 3, 3, 3 \\
9b &$A_7$&$(7.1),\, (1,1,1)$&$2\frac{1}{4}(1,1)$, $A_7$& 1, 1, 8 \\
9c &$A_2$, $A_5$&$(7.3),\, (1,1,1)$&$\frac{1}{9}(1,2)$, $A_2$, $A_5$& 1, 3, 6 \\
9d &$A_1$, $2A_3$&$(7.2),\, (1,1,1)$&$\frac{1}{8}(1,3)$, $2A_3$& 2, 4, 4 \\
10b &$A_8$&$(8.1),\, (1,1,1)$&$2\frac{1}{9}(1,2)$, $A_8$& 1, 1, 9 \\
10c &$A_1$, $A_7$&$(8.2),\, (1,1,1)$&$\frac{1}{16}(1,3)$, $\frac{1}{8}(1,3)$, $A_7$ & 1, 2, 8 \\
10d &$A_1$, $A_2$, $A_5$&$(8.3),\, (1,1,1)$&$\frac{1}{18}(1,5)$, $\frac{1}{12}(1,5)$, $A_5$& 2, 3, 6 \\
10e &$A_4$, $A_4$&$(8.4),\, (1,2,1)$&$\frac{1}{25}(1,9)$, $\frac{1}{20}(1,9)$, $A_4$ & 1, 5, 5 \\ 
11 &$\frac{1}{4}(1,1)$&$1,\, (1,1,2)$&$\frac{1}{4}(1,1)$& 1, 1, 1 \\
15 &$\frac{1}{4}(1,1)$, $A_4$&$4,\, (1,2,1)$&$\frac{1}{4}(1,1)$, $A_4$& 1, 1, 5 \\
18 &$\frac{1}{4}(1,1)$, $A_1$, $A_5$&$6.2,\, (1,1,1)$&$\frac{1}{4}(1,1)$, $A_1$, $A_5$& 1, 2, 6 \\
19 &$\frac{1}{4}(1,1)$, $A_7$&$7.1,\, (1,1,1)$&$2\frac{1}{4}(1,1)$, $A_7$& 1, 1, 8 \\
21 &$\frac{1}{8}(1,3)$, $A_2$&$3,\, (1,2,1)$&$\frac{1}{8}(1,3)$, $A_2$& 1, 2, 3 \\
25 &$2\frac{1}{4}(1,1)$, $A_7$&$7.1,\,(1,1,1)$&$2\frac{1}{4}(1,1)$, $A_7$& 1, 1, 8 \\
26 &$\frac{1}{8}(1,3)$, $2A_3,$&$7.2,\,(1,1,1)$&$\frac{1}{8}(1,3)$, $2A_3$& 2, 4, 4 \\
27 &$\frac{1}{8}(1,3)$, $A_7$&$8.2,\,(1,1,1)$&$\frac{1}{16}(1,3)$, $\frac{1}{8}(1,3)$, $A_7$& 1, 2, 8 \\
30 &$\frac{1}{12}(1,5)$, $2A_2$&$6.1,\, (1,1,2)$&$\frac{1}{12}(1,5)$, $2A_2$& 3, 3, 3 \\
33 &$A_1$, $\frac{1}{12}(1,5)$, $A_5$&$8.3,\, (1,1,1)$&$\frac{1}{18}(1,5)$, $\frac{1}{12}(1,5)$, $A_5$& 2, 3, 6 \\
40 &$\frac{1}{20}(1,9)$&$4, (1,3,2)$&$\frac{1}{9}(1,2)$, $\frac{1}{20}(1,9)$& 1, 1, 5 \\
44 &$\frac{1}{20}(1,9)$, $A_4$&$8.4,\, (1,2,1)$&$\frac{1}{25}(1,9)$, $\frac{1}{20}(1,9)$, $A_4$& 1, 5, 5 \\
46 &$A_2$, $\frac{1}{24}(1,11)$&$7.3,\, (1,2,1)$&$\frac{1}{9}(1,2)$, $A_2$, $\frac{1}{24}(1,11),$& 1, 3, 6 \\
50 &$\frac{1}{36}(1,17)$&$8.1,\, (2,1,1)$&$2\frac{1}{9}(1,2)$, $\frac{1}{36}(1,17)$& 1, 1, 9 \\
\end{longtable}

\section{Existence of special fibrations}

Let $X$ be a log del Pezzo surface such that $\rho(X)=1$ and let $\pi \colon \tilde{X} \rightarrow X$ be its minimal resolution.
We show that, under certain hypotheses, $\tilde{X}$ admits a morphism $p \colon \tilde{X} \rightarrow \bP^1$ with
general fibre a smooth rational curve such that the exceptional locus of $\pi$ has a particularly simple form with respect to
the ruling $p$. 
When $X$ has only $T$-singularities (and satisfies the hypotheses), we use this structure to construct a toric surface $Y$ 
such that $X$ is a $\bQ$-Gorenstein deformation of $Y$, see Sec.~\ref{pfmainthm}.

We first establish the existence of a so called $1$-complement of $K_X$.
We recall the definition and basic properties. For more details and motivation, see \cite[Sec.~19]{FA}, \cite{P}.
Let $X$ be a projective surface with quotient singularities.
A \emph{$1$-complement of $K_X$} is a divisor $D \in |-K_X|$ such that the pair $(X,D)$ is log canonical.
In particular, by the classification of log canonical singularities of pairs \cite[Thm.~4.15]{KM}, $D$ is a nodal curve,
and, at each singularity $P \in X$, either $D=0$ and $P \in X$ is a Du Val singularity, 
or the pair $(P \in X, D)$ is locally analytically isomorphic to the pair $(\frac{1}{n}(1,a),(uv=0))$ for some $n$ and $a$.
Moreover $D$ has arithmetic genus $1$ because $2p_a(D)-2=(K_X+D) \cdot D = 0$ (note that the adjunction formula holds
because $K_X+D$ is Cartier \cite[16.4.3]{FA}). 
Thus $D$ is either a smooth elliptic curve or a cycle of smooth rational curves.

\begin{thm} \label{1complement}
Let $X$ be a log del Pezzo surface such that $\rho(X)=1$.
Assume that $\dim |-K_X| \ge 1$ and every singularity of $X$ is either a cyclic quotient singularity or a Du~Val singularity.
Then there exists a $1$-complement of $K_X$, i.e., a divisor $D \in |-K_X|$ such that the pair $(X,D)$ is log canonical.
\end{thm}
\begin{proof}
Write $-K_X \sim M+F$ where $M$ is an irreducible divisor such that $\dim |M| > 0$ and $F$ is effective 
(we do not assume that $F$ is the fixed part of $|-K_X|$).
Let $M$ be general in $|M|$.

Suppose first that $(X,M)$ is purely log terminal (plt).
Then $M$ is a smooth curve.
We may assume that $F \neq 0$ (otherwise $M$ is a $1$-complement).
Then $-(K_X+M) \sim F $ is ample (because $\rho(X)=1$).
Recall that for $X$ a normal variety and $S \subset X$ an irreducible divisor the \emph{different} $\Diff_S(0)$ is the effective $\bQ$-divisor on $S$ defined by the equation
$$(K_X+S)|_S = K_S + \Diff_S(0).$$
That is, $\Diff_S(0)$ is the correction to the adjunction formula for $S \subset X$ due to the singularities of $X$ at $S$.
See \cite[Sec.~16]{FA}.
If $S$ is a normal variety and $B$ is an effective $\bQ$-divisor on $S$ with coefficients less than $1$, a \emph{1-complement}
of $K_S+B$ is a divisor $D \in |-K_S|$ such that $(S,D)$ is log canonical and $D \ge \lfloor 2B \rfloor$.
By \cite[Prop.~4.4.1]{P} it's enough to show that $K_M+\Diff_M(0)$ has a $1$-complement.

The curve $M$ is smooth and rational and $\deg(K_M+\Diff_M(0))<0$ because
$$
2p_a(M)-2 \le \deg(K_M+\Diff_M(0)) = (K_X+M)\cdot M = -F \cdot M < 0.
$$
Moreover, at each singular point $P_i$ of $X$ on $M$, the pair $(X,M)$ is of the form $(\frac{1}{m_i}(1,a_i),(x=0))$, 
and 
$$\Diff_M(0)=\sum_i \left(1-\frac{1}{m_i}\right)P_i$$
by \cite[16.6.3]{FA}. 
So, if $K_M+\Diff_M(0)$ does not have a $1$-complement, then, by \cite[Cor.~19.5]{FA} or direct calculation,
there are exactly $3$ singular points of $X$ on $M$, 
and $(m_1,m_2,m_3)$ is a Platonic triple $(2,2,m)$ (for some $m \ge 2$), $(2,3,3)$, $(2,3,4)$, or $(2,3,5)$. 
The divisor $F$ passes through each singular point $P_i$ because $F \sim -(K_X+M)$ is not Cartier there.
So $F\cdot M \ge \sum \frac{1}{m_i}$, 
and 
$$0=(K_X+M+F)\cdot M = \deg(K_M+\Diff_M(0))+ F\cdot M  \ge 1,$$
a contradiction.

Now suppose that the pair $(X,M)$ is not plt, and let $c$ be its log canonical threshold, i.e.,
$$c = \sup \, \{ t \in \bQ_{\ge 0} \ | \ (X,tM) \mbox{ is log canonical } \}.$$
Then there exists a projective birational morphism $f \colon Y \rightarrow X$ with exceptional locus an irreducible divisor $E$ such
that the discrepancy $a(E,X,cM)=-1$ and $(Y,E)$ is plt. See \cite[Prop.~3.1.4]{P}. 
So 
$$K_Y+cM'+E=f^*(K_X+cM)$$ 
where $M'$ is the strict transform of $M$.
Now 
$$-(K_Y+E)=cM'-f^*(K_X+cM)$$ 
is nef (note $M'$ is nef because it moves).
Moreover $-(K_Y+E)$ is big unless $M'^2=0$ and $K_X+cM \sim_{\bQ} 0$, in which case $c=1$, $F=0$, and $M$ is a $1$-complement.
So we may assume $-(K_Y+E)$ is nef and big.
Thus, by \cite[Prop.~4.4.1]{P} again, it's enough to show that $K_E+\Diff_E(0)$ has a $1$-complement.
Suppose not. 
Then $E$ passes through $3$ cyclic quotient singularities on $Y$ as above.
Let $\tilde{Y} \rightarrow Y$ be the minimal resolution of $Y$, $E'$ the strict transform of $E$, and consider the composition
$g \colon \tilde{Y} \rightarrow X$.
Let $P \in X$ be the point $f(E)$.
Then $g^{-1}(P)$ is the union of $E'$ and $3$ chains of smooth rational curves (the exceptional loci of the minimal resolutions of the cyclic 
quotient singularities), and $E'$ meets each chain in one of the end components.
Let $-b_i$ be the self-intersection number of the end component $F_i$ of the $i$th chain that meets $E'$.
Then $b_i \le m_i$ where $m_i$ is the order of the cyclic group for the $i$th quotient singularity.
If we contract the $F_i$ and let $\oE'$ denote the image of $E'$, then 
$$0 > {\oE'}^2={E'}^2 + \sum \frac{1}{b_i} \ge {E'}^2 + \sum \frac{1}{m_i} > {E'}^2+1.$$
Hence $E'^2 \le -2$ and $g$ is the minimal resolution of $X$.
So $P \in X$ is a $D$ or $E$ singularity by our assumption. 
But $P \in X$ is a basepoint of $|-K_X|$, so this contradicts Thm.~\ref{DorE}.
\end{proof}

We describe the types of degenerate fibres which occur in the ruling we construct.
We first introduce some notation.

\begin{defn} 
Let $a,n \in \bN$ with $a<n$ and $(a,n)=1$. We say the fractions $n/a$ and $n/(n-a)$ are \emph{conjugate}. 
\end{defn}

\begin{lem} \label{conjugate} 
If $[b_1,\ldots,b_r]$ and $[c_1,\ldots,c_s]$ are conjugate, then so are $[b_1+1,b_2,\ldots,b_r]$
and $[2,c_1,\ldots,c_s]$. Conversely, every conjugate pair can be constructed from $[2]$,$[2]$ by a sequence of such 
steps. Also, if $[b_1,\ldots,b_r]$ and $[c_1,\ldots,c_s]$ are conjugate then so are $[b_r,\ldots,b_1]$ and $[c_s,\ldots,c_1]$.
\end{lem}
\begin{proof}
If $[b_1,\ldots,b_r]=n/a$ and $[c_1,\ldots,c_s]=n/(n-a)$ then $[b_1+1,b_2,\ldots,b_r]=(n+a)/a$ and
$[2,c_1,\ldots,c_s]=(n+a)/n$. The last statement follows immediately from Rem.~\ref{reverse}.
\end{proof}

\begin{prop} \label{fibres}
Let $S$ be a smooth surface, $T$ a smooth curve, and $p \colon S \rightarrow T$ a morphism with general fibre a smooth rational curve.
Let $f$ be a degenerate fibre of $p$. 
Suppose that $f$ contains a unique $(-1)$-curve and the union of the remaining irreducible components of $f$ 
is a disjoint union of chains of smooth rational curves. Then the dual graph of $f$ has one of the following forms.

\[
\xymatrix
@R=1pc
@C=1.8pc
{
\overset{a_r}\circ\ar@{-}[r]&\cdots\ar@{-}[r]&
\overset{a_1}\circ\ar@{-}[r]&\bullet\ar@{-}[r]
&\overset{b_1}\circ\ar@{-}[r]&\cdots\ar@{-}[r]&
\overset{b_s}\circ
}
\leqno{(I)}
\]

\[
\xymatrix
@R=1pc
@C=1.8pc
{
\overset{a_r}\circ\ar@{-}[r]&\cdots\ar@{-}[r]&
\overset{a_1}\circ\ar@{-}[r]&\overset{t+2}\circ\ar@{-}[r]
&\overset{b_1}\circ\ar@{-}[r]&\cdots\ar@{-}[r]&
\overset{b_s}\circ
\\
&&&\bullet\ar@{-}[u]\ar@{-}[r]&\underset{2}{\circ}
\ar@{-}[r]&\cdots\ar@{-}[r]
&\underset{2}{\circ}&
}
\leqno{(II)}
\]
Here the black vertex denotes the $(-1)$-curve and a white vertex with label $a \ge 2$
denotes a smooth rational curve with self-intersection number $-a$.
In both types the strings $[a_1,\ldots,a_r]$ and $[b_1,\ldots,b_s]$ are conjugate.
In type $(II)$ there are $t$ $(-2)$-curves in the branch containing the $(-1)$-curve.

Conversely, any configuration of curves of this form is a degenerate fibre of a fibration $p \colon S \rightarrow T$ as above.
\end{prop}

\begin{proof}
The morphism $p \colon S \rightarrow T$ is obtained from a $\bP^1$-bundle $F \rightarrow T$ by a sequence of blowups.
The statements follow by induction on the number of blowups.
\end{proof}

We refer to the fibres above as fibres of types $(I)$ and $(II)$. We also call a fibre of the form
\[
\xymatrix{
\bullet\ar@{-}[r]&\overset{2}\circ\ar@{-}[r]&\cdots\ar@{-}[r]&\overset{2}\circ\ar@{-}[r]&
\bullet
}
\leqno{(O)}
\]
a fibre of type $(O)$.

\begin{rem} \label{remfibres2}
The curves of multiplicity one in the fibre are the ends of the chain in types $(O)$ and $(I)$
and the ends of the branches not containing the $(-1)$-curve in type $(II)$.
In particular, a section of the fibration meets the fibre in one of these curves.
\end{rem} 

\begin{thm} \label{fibration}
Let $X$ be a log del Pezzo surface such that $\rho(X)=1$.
Assume that $\dim |-K_X| \ge 1$ and every singularity of $X$ is either a cyclic quotient singularity or a Du~Val singularity.
Let $\pi \colon \tilde{X} \rightarrow X$ be the minimal resolution of $X$.
Then one of the following holds.
\begin{enumerate}
\item There exists a morphism $p \colon \tilde{X} \rightarrow \bP^1$ with general fibre a smooth rational curve satisfying 
one of the following.
\begin{enumerate}
\item Exactly one component $\tilde{E}_1$ of the exceptional locus of $\pi$ is $p$-horizontal. The curve 
$\tilde{E_1}$ is a section of $p$. The fibration $p$ has at most two degenerate fibres and each is of type $(I)$ or $(II)$.
\item Exactly two components $\tilde{E}_1,\tilde{E_2}$ of the exceptional locus of $\pi$ are $p$-horizontal.
The curves $\tilde{E_1}, \tilde{E_2}$ are sections of $p$. 
Either $\tilde{E_1}$ and $\tilde{E_2}$ are disjoint and $p$ has two degenerate fibres of types $(O)$ and either $(I)$ or $(II)$,
or $\tilde{E_1} \cdot \tilde{E_2}=1$ and $p$ has a single degenerate fibre of type $(O)$.
The sections $\tilde{E_1}$ and $\tilde{E_2}$ meet distinct components of the degenerate fibres. 
\end{enumerate}
\item The surface $X$ has at most $2$ non Du~Val singularities and each is of the form $\frac{1}{m}(1,1)$ for some $m \ge 3$.
\end{enumerate}
\end{thm}
\begin{proof}
Assume that $K_X$ is not Cartier.
As in the proof of Thm.~\ref{DorE}, let $\nu \colon \hat{X} \rightarrow X$ be the minimal resolution of the non Du Val singularities,
$\{ E_i \}$ the exceptional divisors, and $E=\sum E_i$.
Write $|-K_{\hat{X}}| =|M|+F$ where $F$ is the fixed part and $M \in |M|$ is general.
Then $F \ge E$ and $\dim|M|=\dim|-K_X| \ge 1$.

We run the MMP on $\hat{X}$.
We obtain a birational morphism $\phi \colon \hat{X} \rightarrow \oX$ such that $\oX$ has Du Val singularities and 
either $\rho(\oX)=2$ and there is a fibration
$\psi \colon \oX \rightarrow \bP^1$ such that $-K_{\oX}$ is $\psi$-ample or $\rho(\oX)=1$ and $-K_{\oX}$ is ample.
Moreover, $\phi$ is a composition
$$\hat{X}=X_1 \stackrel{\phi_1}{\longrightarrow} X_2 \stackrel{\phi_2}{\longrightarrow} \ldots \stackrel{\phi_n}{\longrightarrow}X_{n+1}=\oX$$
where $\phi_i$ is a weighted blowup of a smooth point of $X_{i+1}$ with weights $(1,n_i)$ 
(by the classification of birational extremal contractions in the MMP for surfaces with Du Val singularities).
\begin{claim}
Given $\phi \colon \hat{X} \rightarrow \oX$, we can direct the MMP so that the components of $E$ contracted by $\phi$ are contracted last. That is, for some $1 \le m \le n$, the exceptional divisor of $\phi_i$ is (the image of) a component of $E$ iff $i>m$.
\end{claim}
\begin{proof}
We have $K_{\hat{X}}=\nu^*K_X+\sum a_iE_i$ where $-1<a_i<0$ for each $i$.
Write $\Delta= \sum (-a_i)E_i$.
So $\nu^*K_X=K_{\hat{X}}+\Delta$ and $\Delta$ is an effective divisor such that $\lfloor \Delta \rfloor=0$ and $\Supp \Delta =E$.
Hence $-(K_{\hat{X}}+\Delta)$ is nef and big and $(\hat{X},\Delta)$ is Kawamata log terminal (klt).
These properties are preserved under the $K_{\hat{X}}$-MMP.

Let $R=\sum R_i$ be the sum of the $\phi$-exceptional curves that are not contained in $E$ and $R' \subset R$ a connected component.
Then $R' \cdot E > 0$ (otherwise $\nu$ is an isomorphism near $R'$ which contradicts $\rho(X)=1$).
Let $R_i$ be a component of $R'$ such that $R_i \cdot E > 0$.
Then $(K_{\hat{X}}+\Delta)\cdot R_i \le 0$ and $R_i \cdot \Delta > 0$. 
So $K_{\hat{X}} \cdot R_i < 0$, and we can contract $R_i$ first in the $K_{\hat{X}}$-MMP.
Repeating this procedure, we contract all of $R$, obtaining a birational morphism $\hat{X} \rightarrow \hat{X}'$.
Finally we run the MMP on $\hat{X}'$ over $\oX$ to contract the remaining curves.
\end{proof}
\begin{claim}
We may assume $\rho(\oX)=2$.
\end{claim}
\begin{proof}
Suppose $\rho(\oX)=1$.
Write $\oM =\phi_*M$, etc.
Then $-K_{\oX} \sim \oM + \oF$, $\oF \ge \oE > 0$, and $\oE$ is contained in the smooth locus of $\oX$.
Thus, as in the proof of Thm.~\ref{DorE}, $-K_{\oX}$ is not a generator of $\Pic \oX$, 
so $\oX \simeq \bP^2$ or $\oX \simeq \bP(1,1,2)$ by the classification of log del Pezzo surfaces with Du Val singularities.
In particular, it follows that $\oE$ has at most $2$ components. 

Suppose first that $\phi$ does not contract any component of $E$. Then $E$ has at most $2$ components.
So, either we are in case (2), or $E=E_1+E_2$, $E_1 \cap E_2 \neq \emptyset$, $\oX \simeq \bP^2$,
and $\oM, \oE_1, \oE_2 \sim l$, where $l$ is the class of a line.
In this case $\rho(\hat{X})= \rho(X)+2=3$, so $\phi \colon \hat{X} \rightarrow \oX$ is a composition of two weighted blowups of weights 
$(1,n_1)$, $(1,n_2)$.
These must have centres two distinct points $P_1 \in \oE_1$, $P_2 \in \oE_2$, 
and in each case the local equation of $\oE_i$ is a coordinate with weight $n_i$ (because $E_i$ is contained in the smooth locus of $\hat{X}$).
Let $l_{12}$ be the line through $P_1$ and $P_2$.
Then these blowups are toric with respect to the torus $\oX \setminus l_{12}+\oE_1+\oE_2$.
We find that the minimal resolution $\tilde{X}$ is a toric surface with boundary divisor a cycle of smooth rational curves with 
self-intersection numbers 
$$-2,\ldots,-2,-1,-(n_1-1),-(n_2-1),-1,-2,\ldots,-2,-1$$ 
where $\tilde{E}_1$ and $\tilde{E}_2$ are the curves with self-intersection numbers $-(n_1-1),-(n_2-1)$,
the first two $(-1)$-curves are the strict transforms of the exceptional curves of the blowups of $P_1$ and $P_2$, 
the last $(-1)$-curve is the strict transform of $l_{12}$,
and the chains of $(-2)$-curves are the exceptional loci of the resolutions of the singularities of $\hat{X}$ and have lengths $(n_1-1)$ and $(n_2-1)$.
In particular, there is a fibration $p \colon \tilde{X} \rightarrow \bP^1$ with two degenerate fibres of types 
$-1,-2,\ldots,-2,-1$ (where there are $(n_2-1)$ $(-2)$-curves) and $-2,\ldots,-2,-1,-(n_1-1)$ (where there are $(n_1-2)$ $(-2)$-curves),
and two $\pi$-exceptional sections with self-intersection numbers $-(n_2-1)$ and $-2$. So we are in case~(1b).

Now suppose $\phi$ contracts some component of $E$.
Then $\phi_n \colon X_n \rightarrow X_{n+1}=\oX$ is an (ordinary) blowup of a smooth point $Q \in \oX$.
If $\oX \simeq \bP^2$ then $X_n \simeq \bF_1$ and there is a fibration $\psi \colon X_n \rightarrow \bP^1$.
So we may assume $\rho(\oX)=2$.
If $\oX \simeq \bP(1,1,2)$, the quadric cone, let $L$ be the ruling of the cone through $Q$.
Then the strict transform $L'$ of $L$ on $X_n$ satisfies $K_{X_n} \cdot L' < 0$ and ${L'}^2 <0$.
Contracting $L'$ we obtain a morphism $\phi_n' \colon X_n \rightarrow \oX' \simeq \bP^2$.
So, replacing $\phi_n$ by $\phi_n'$, we may assume $\oX \simeq \bP^2$.
\end{proof}

We now assume $\rho(\oX)=2$.
We have a diagram
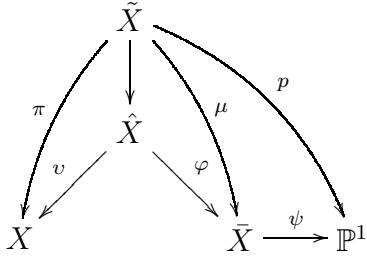
\[
\xymatrix{
&\tilde X\ar@/^01.2pc/[rrdd]^{p}\ar@/^0.7pc/[rdd]^{\mu}
\ar@/_0.7pc/[ldd]_{\pi}\ar[d]^{}&
\\
&\hat X\ar[ld]_{\upsilon}\ar[rd]^{\varphi}&
\\
X&&\bar X\ar[r]^{\psi}&\bP^1
}
\]
where $\pi \colon \tilde{X} \rightarrow X$ is the minimal resolution.
Let $l$ be a general fibre of $p$ and $\tilde{E}$ the strict transform of $E$ on $\tilde{X}$.
Note that, by construction, the components of the exceptional locus of $\pi$ over Du Val singularities are contained in fibres of $p$.
Write $|-K_{\tilde{X}}|=|\tilde{M}|+\tilde{F}$ where $\tilde{F}$ is the fixed part and $\tilde{M} \in |\tilde{M}|$ is general.
Then $\tilde{F} \ge \tilde{E}$.

There is a $1$-complement of $K_X$ by Thm.~\ref{1complement}.
This can be lifted to $\tilde{X}$.
(Indeed, if $D$ is a $1$-complement of $K_X$, define $\tilde{D}$ by $K_{\tilde{X}}+\tilde{D}=\pi^*(K_X+D)$ and $\pi_*\tilde{D}=D$.
Note that $\tilde{D}$ is an effective $\bZ$-divisor because $K_{\tilde{X}}$ is $\pi$-nef and $K_X+D$ is Cartier.
Then $\tilde{D}$ is a $1$-complement of $K_{\tilde{X}}$.)
Hence $(\tilde{X},\tilde{M}+\tilde{F})$ is log canonical. In particular, $\tilde{F}$ is reduced and $\tilde{M}+\tilde{F}$ is a cycle of smooth rational curves.

There exists a $p$-horizontal component $\tilde{E_1}$ of $\tilde{E}$ (because $\rho(X)=1$).
Then 
$$1 \le \tilde{E_1}\cdot l \le (\tilde{F}+\tilde{M}) \cdot l = -K_{\tilde{X}} \cdot l =2.$$

Suppose first that $\tilde{E_1} \cdot l=2$.
Then $\tilde{M}$ and $\tilde{F}-\tilde{E_1}$ are $p$-vertical.
Hence $\tilde{M} \sim l$ and $\tilde{F}=\tilde{E}_1$, so $\tilde{E}=\tilde{E}_1$ and we are in case~(2).

Suppose now that $\tilde{E_1} \cdot l=1$.
Since $\mu(\tilde{E}_1)$ is contained in the smooth locus of $\oX$, the fibres of $\psi$ have multiplicity $1$, so $\psi$ is smooth
by \cite[Lem.~11.5.2]{KMcK}.
Thus $\oX \simeq \bF_n$ for some $n \ge 0$.

If $\tilde{E}_1$ is the only $p$-horizontal component of $\tilde{E}$ we are in case~(1a).
Suppose there is another $p$-horizontal component $\tilde{E}_2$.
Then, since $-K_{\tilde{X}}\cdot l =2$, we have $\tilde{E_2} \cdot l =1$ and $\tilde{M}$ and $\tilde{F}-\tilde{E}_1-\tilde{E}_2$ are contained in fibres of $p$.
If $M \sim 2l$ then $\tilde{F}=\tilde{E}=\tilde{E}_1+\tilde{E}_2$ and $\tilde{E_1} \cap \tilde{E_2}= \emptyset$ so we are in case~(1b).
So we may assume $\tilde{M} \sim l$.
Then the components of $\tilde{F}$ form a chain, with ends $\tilde{E}_1$ and $\tilde{E_2}$.

We note that a component $\Gamma$ of a degenerate fibre of $p$ that is not contracted by $\pi$ is necessarily a $(-1)$-curve,
because $K_{\tilde{X}}=\pi^*K_X - \tilde{\Delta}$ where $\tilde{\Delta}$ is effective and $\pi$-exceptional, so
$$K_{\tilde{X}} \cdot \Gamma \le \pi^*K_X \cdot \Gamma = K_X \cdot \pi_* \Gamma < 0.$$
Hence, since $\rho(X)=1$, there exists a unique degenerate fibre of $p$ containing exactly two $(-1)$-curves, and any other degenerate
fibres contain exactly one $(-1)$-curve.
Let $\tilde{G}$ denote the reduction of the fibre containing two $(-1)$-curves.

If $\tilde{F}=\tilde{E}_1+\tilde{E}_2$ then $\tilde{E_1} \cdot \tilde{E_2}=1$ and 
any degenerate fibre of $p$ consists of $(-1)$-curves and $(-2)$-curves.
It follows that $\tilde{G}$ is of type $(O)$ and there are no other degenerate fibres, so we are in case~(1b).
So assume $\tilde{F} > \tilde{E}_1+\tilde{E}_2$. Then $\tilde{E}_1 \cap \tilde{E_2} = \emptyset$.

Suppose first that $\tilde{G}$ is the only degenerate fibre. Then $\tilde{F} \le \tilde{G} +\tilde{E}_1+\tilde{E}_2$.
Write $\tilde{G}=\tilde{G}'+\tilde{G}''$ where $\tilde{G}'=\tilde{F}-\tilde{E}_1-\tilde{E}_2$.
So $\tilde{G}'$ is a chain of smooth rational curves.
It follows that each connected component of $\tilde{G}''$ is a chain of smooth rational rational curves such that one end component 
is a $(-1)$-curve adjacent to $\tilde{G}'$ and the remaining curves are $(-2)$-curves.
We construct an alternative ruling $p' \colon \tilde{X} \rightarrow \bP^1$ with only one horizontal $\pi$-exceptional curve by inductively contracting $(-1)$-curves as follows.
First contract the components of $\tilde{G}''$.
Second, contract $(-1)$-curves in $\tilde{G}'$ until the image of $\tilde{E}_1$ or $\tilde{E}_2$ is a $(-1)$-curve.
Now contract this curve, and continue contracting $(-1)$-curves until we obtain a ruled surface $\oX' \simeq \bF_m$.
Then $\tilde{M} \sim l$ is horizontal for the induced ruling $p'$.
Moreover, if $C$ is a $p'$-horizontal $\pi$-exceptional curve then $C \not \subset \tilde{G}''$ by construction.
Hence $C \subset \tilde{F}$. Thus there exists a unique such $C$, and $C$ is a section of $p'$.
So we are in case~(1a).

Finally, suppose there is another degenerate fibre of $p$, and let $\tilde{V}$ denote its reduction.
Then $\tilde{V}$ contains a unique $(-1)$-curve $C$.
The surface $X$ has only cyclic quotient singularities by assumption.
Therefore $\tilde{V}-C$ is a union of chains of smooth rational curves.
It follows that $\tilde{V}$ is a fibre of type $(I)$ or $(II)$.
Now $\tilde{E_1} \cdot C = \tilde{E_2} \cdot C = 0$ because $C$ has multiplicity greater than $1$ in the fibre.
So $\tilde{V}$ contains a component of $\tilde{F}$ (because $1=-K_{\tilde{X}} \cdot C =(\tilde{M}+\tilde{F}) \cdot C$).
Hence $\tilde{F}-\tilde{E}_1-\tilde{E}_2 \le \tilde{V}$ (because $\tilde{M}+\tilde{F}$ is a cycle of rational curves and 
$\tilde{M} \sim l$).
In particular, $\tilde{G}$ consists of two $(-1)$-curves and some $(-2)$-curves.
Hence $\tilde{G}$ is of type $(O)$ and we are in case~(1b). This completes the proof.
\end{proof}

\section{Proof of Main Theorem}\label{pfmainthm}

\begin{thm} \label{deformation}
Let $X$ be a log del Pezzo surface such that $\rho(X)=1$ and $X$ has only $T$-singularities.
Then exactly one of the following holds
\begin{enumerate}
\item $X$ is a $\bQ$-Gorenstein deformation of a toric surface $Y$, or
\item $X$ is one of the sporadic surfaces listed in Example~\ref{sporadic}.
\end{enumerate}
\end{thm}

\begin{rem}
Note that the surface $Y$ in Thm.~\ref{deformation}$(1)$ necessarily has only $T$-singularities and $\rho(Y)=1$.
Thus $Y$ is one of the surfaces listed in Thm.~\ref{toric}.
\end{rem}

\begin{ex} \label{sporadic}
We list the log del Pezzo surfaces $X$ such that $X$ has only $T$-singularities and $\rho(X)=1$, but $X$ is not a $\bQ$-Gorenstein 
deformation of a toric surface. In each case $X$ has index $\le 2$.
If $X$ is Gorenstein, the possible configurations of singularities are 
$$D_5,\, E_6,\, E_7,\, A_1D_6,\, 3A_1D_4,\, E_8,\, D_8,\, A_1E_7,$$
$$ A_2E_6,\, 2A_1D_6, \, A_3D_5,\, 2D_4,\, 2A_1 2A_3,\, 4A_2.$$
The configuration determines the surface uniquely with the following exceptions: there are two surfaces for $E_8$, $A_1E_7$, $A_2E_6$,
and an $\bA^1$ of surfaces for $2D_4$. See \cite[Thm~4.3]{AN}. 
If $X$ has index $2$, the possible configurations of singularities are 
$$\frac{1}{4}(1,1)D_8, \frac{1}{4}(1,1)2A_1D_6, \frac{1}{4}(1,1)A_3D_5, \frac{1}{4}(1,1)2D_4,$$
and the configuration determines the surface uniquely. See \cite[Thm~4.2]{AN}.
\end{ex}

\begin{rem}
Note that the case $K_X^2=7$ does not occur. This may be explained as follows.
If $X$ is a del Pezzo surface with $T$-singularities such that $\rho(X)=1$,
then there exists a $\bQ$-Gorenstein smoothing $\cX/T$ of $X$ over $T:= \Spec k[[t]]$ such that the generic fibre 
$\cX_K$ is a smooth  del Pezzo surface over $K=k((t))$ with $\rho(\cX_K)=1$.
(Indeed, if $\cX/T$ is a smoothing of $X$ over $T$,
the restriction map $\Cl(\cX) \rightarrow \Cl(\cX_K)=\Pic(\cX_K)$ is an isomorphism because the closed fibre $X$ is irreducible
and the restriction map $\Pic(\cX) \rightarrow \Pic(X)$ is an isomorphism because $H^1(\cO_X)=H^2(\cO_X)=0$.
Thus $\rho(\cX_K) \ge \rho(X)=1$ with equality iff the total space $\cX$ of the deformation is $\bQ$-factorial.
Since there are no local-to-global obstructions for deformations of $X$, there exists a $\bQ$-Gorenstein smoothing $\cX/T$ such that
$P \in \cX$ is smooth for $P \in X$ a Du Val singularity and 
$P \in \cX$ is of type $\frac{1}{n}(1,-1,a)$ for $P \in X$ a singularity of type $\frac{1}{dn^2}(1,dna-1)$ (see Sec.~\ref{Tdefclass}). 
In particular, $\cX$ is $\bQ$-factorial.)  
Note that $K_{\cX_{K}}^2=K_X^2$. If $Y$ is a smooth del Pezzo surface with $K_Y^2=7$ over a field (not necessarily algebraically closed) then $\rho(Y)>1$, see, e.g., \cite{Ma}. Hence there is no $X$ with $K_X^2=7$.
\end{rem} 

\begin{proof}[Proof of Thm.~\ref{deformation}]
First assume that $X$ does not have a $D$ or $E$ singularity. 
Note that $\dim |-K_X| = K_X^2 \ge 1$ by Lem.~\ref{RR}, so we may apply Thm.~\ref{fibration}. 
We use the notation of that theorem.

Suppose first that we are in case $(1a)$. We construct a toric surface $Y$ and prove that $X$ is a $\bQ$-Gorenstein deformation of $Y$.
We first describe the surface $Y$. Let $\tilde{E}_1^2=-d$.
There is a uniquely determined toric blowup $\mu_Y \colon \tilde{Y} \rightarrow \bF_d$ such that
$\mu_Y$ is an isomorphism over the negative section $B \subset \bF_d$, and the degenerate fibres of the ruling 
$p_Y \colon \tilde{Y} \rightarrow \bP^1$ are fibres of type $(I)$ associated to the degenerate fibres of 
$p \colon \tilde{X} \rightarrow \bP^1$ as follows.
Let $f$ be a degenerate fibre of $p$ of type $(I)$ or $(II)$ as in Prop.~\ref{fibres}, 
and assume that $\tilde{E}_1$ intersects the left end component.
If $f$ is of type $(I)$ then the associated fibre $f_Y$ of $p_Y$ has the same form.
If $f$ is of type $(II)$ then $f_Y$ is a fibre of type $(I)$ with self-intersection numbers 
$$-a_r,\ldots,-a_1,-t-2,-b_1,\ldots,-b_s,-1,-d_1,\ldots,-d_u$$ 
Note that the sequence $d_1,\ldots,d_u$ is uniquely determined (see Prop.~\ref{fibres}).
In each case the strict transform $B'$ of $B$ again intersects the left end component of $f_Y$.

Let $Y$ be the toric surface obtained from $\tilde{Y}$ by contracting the strict transform of the negative section of $\bF_d$ 
and the components of the degenerate fibres of the ruling with self-intersection number at most $-2$.
For each fibre $f$ of $p$ of type $(II)$ as above, the chain of rational curves with self-intersections $-d_1,\ldots,-d_u$
in the associated fibre $f_Y$ of $p_Y$  contracts to a $T_{t+1}$ singularity by Lem.~\ref{T}(1).
This singularity replaces the $A_t$ singularity on $X$ obtained by contracting the chain of $t$ $(-2)$-curves in $f$. 
In particular, the surface $Y$ has $T$-singularities. Moreover $\rho(Y)=1$, and $K_Y^2=K_X^2$ by Prop.~\ref{noether}. 
A $T_d$-singularity admits a $\bQ$-Gorenstein deformation to an $A_{d-1}$ singularity (see Prop.~\ref{Tdeformation}). 
Hence the singularities of $X$ are a $\bQ$-Gorenstein deformation of the singularities of $Y$.
There are no local-to-global obstructions for deformations of $Y$ by Prop.~\ref{unobs}. 
Hence there is a $\bQ$-Gorenstein deformation $X'$ of $Y$ with the same singularities as $X$.
We prove below that $X \simeq X'$.
 
Let $f$ be a degenerate fibre of $p$ of type $(II)$ as above and $f_Y$ the associated fibre of $p_Y$.
Let $P \in Y$ be the $T$-singularity obtained by contracting the chain 
of rational curves in $f_Y$ with self-intersections $-d_1,\cdots,-d_u$.
Let $X'$ be the general fibre of a $\bQ$-Gorenstein deformation of $Y$ over the germ of a curve
which deforms $P \in Y$ to an $A_{t}$ singularity and is locally trivial elsewhere. 
Let $\hat{Y} \rightarrow Y$ and $\hat{X}' \rightarrow X'$ be the minimal resolutions of the remaining singularities
(where the deformation is locally trivial).
Thus $\hat{Y}$ has a single $T$-singularity and $\hat{X}'$ a single $A_{t}$ singularity.
The ruling $p_Y \colon \tilde{Y} \rightarrow \bP^1$ descends to a ruling $\hat{Y} \rightarrow \bP^1$; let $A$ be a general fibre of this ruling. 
Then $A$ deforms to a $0$-curve $A'$ in $\hat{X}'$ (because $H^1(\cN_{A/\hat{Y}})=H^1(\cO_A)=0$) 
which defines a ruling $\hat{X}' \rightarrow \bP^1$. 
Let $\tilde{X}' \rightarrow \hat{X}'$ be the minimal resolution of $\hat{X}'$ and consider the induced ruling 
$p_{X'} \colon \tilde{X}' \rightarrow \bP^1$. Note that the exceptional locus of $\hat{Y} \rightarrow  Y$ deforms without change by 
construction. Moreover, the $(-1)$-curve in the remaining degenerate fibre (if any) of $p_Y$ also deforms. 
There is a unique horizontal curve in the exceptional locus of $\pi_{X'} : \tilde{X}' \rightarrow X'$, and $\rho(X')=1$ by 
Prop.~\ref{noether}.
Hence each degenerate fibre of $p_{X'}$ contains a unique $(-1)$-curve, and the remaining components of the fibre 
are in the exceptional locus of $\pi_{X'}$. 
We can now describe the degenerate fibres of $p_{X'}$. 
If $p_Y$ has a degenerate fibre besides $f_Y$, then $p_{X'}$ has a degenerate fibre of the same form.
We claim that there is exactly one additional degenerate fibre of $p_{X'}$, which is of type $(II)$ and has the same form as the fibre 
$f$ of $p$.
Indeed, the union of the remaining degenerate fibres consists of the chain of rational curves with self-intersections 
$-a_r,\ldots,-a_1,-t-2,b_1,\ldots,b_s$ (the deformation of the chain of the same form in $f_Y$), 
the chain of $(-2)$-curves which contracts to the $A_t$ singularity, and some $(-1)$-curves.
The claim follows by the description of degenerate fibres in Prop.~\ref{fibres}.
If there is a second degenerate fibre of $p$ of type $(II)$ we repeat this process.
We obtain a $\bQ$-Gorenstein deformation $X'$ of $Y$ with minimal resolution $\pi_{X'} \colon \tilde{X}' \rightarrow X'$, 
and a ruling 
$p_{X'} \colon \tilde{X}' \rightarrow \bP^1$ such that the exceptional locus of $\pi_{X'}$ 
has the same form with respect to the ruling $p_{X'}$ as that of $\pi$ with respect to $p$.

We claim that $X \simeq X'$. Indeed, there is a toric variety $Z$ and, for each fibre $f_i$ of $p$ of type $(II)$, 
an irreducible toric boundary divisor 
$\Delta_i \subset Z$ and points $P_i,P'_i$ in the torus orbit $O_i \subset \Delta_i$, such that $\tilde{X}$ (respectively $\tilde{X}'$) 
is obtained from $Z$ by successively blowing up the points $P_i$ (respectively $P'_i$)  $t_i+1$ times, where $t_i$ is the length of 
the chain of $(-2)$-curves in $f_i$. It remains to prove that we may assume $P_i=P'_i$ for each $i$.
Let $T$ be the torus acting on $Z$ and $N$ its lattice of 1-parameter subgroups. Let $\Sigma \subset N_{\bR}$ be the fan corresponding to 
$X$ and $v_i \in N$ the minimal generator of the ray in $\Sigma$ corresponding to $\Delta_i$. 
Then $T_i = (N/\langle v_i \rangle) \otimes \bG_m$ is the quotient torus of $T$ which acts faithfully on $\Delta_i$.
Thus, there is an element $t \in T$ taking $P_i$ to $P'_i$ for each $i$ except in the following case: 
there are two fibres of $p$ of type $(II)$,
and $v_1+v_2=0$. In this case, there is a toric ruling $q \colon Z \rightarrow \bP^1$ given by the projection 
$N \rightarrow N/\langle v_1 \rangle$. The toric boundary of $Z$ decomposes into two sections (given by $\Delta_1,\Delta_2$) and two fibres 
of $q$. But one of these fibres (the one containing the image of $\tilde{E}_1 \subset \tilde{X}$) is a chain of rational curves of self-intersections at most $-2$, a contradiction.

Next assume that we are in case $(1b)$.
There is a ruling $p \colon \tilde{X} \rightarrow \bP^1$ with two $\pi$-exceptional sections $\tilde{E_1}$ and $\tilde{E_2}$.
Suppose first that $\tilde{E_1} \cap \tilde{E_2} = \emptyset$. Then there are two degenerate fibres of types $(O)$ and either $(I)$ or $(II)$.
We use the notation of Prop.~\ref{fibres}.
The exceptional locus of $\pi$ consists of the components of the degenerate fibres of self-intersection $\le -2$
and the two disjoint sections $\tilde{E_1}$ and $\tilde{E_2}$ of $p$ which meet the first fibre in the two $(-1)$-curves
and the second fibre in the components labelled $-a_r$ and $-b_s$ respectively. 
If the degenerate fibres are of types $(O)$ and $(I)$ then $X$ is toric.
So we may assume the degenerate fibres are of types $(O)$ and $(II)$.
Set $\tilde{E}_1^2=-a_{r+1}$ and $\tilde{E}_2^2=-b_{s+1}$.
Let $m$ be the number of $(-2)$-curves in the fibre of type $(O)$.
Then $X$ has singularities $A_m$, $A_t$, and the cyclic quotient singularity whose minimal resolution has
exceptional locus the chain of rational curves with self-intersections $-a_{r+1},\ldots,-a_1,-(t+2),-b_1,\ldots,-b_{s+1}$.

The ruling $p \colon \tilde{X} \rightarrow \bP^1$ is obtained from a $\bP^1$-bundle by a sequence of blowups. It follows that 
$m=a_{r+1}+b_{s+1}-2$.

We construct a toric surface $Y$ and prove that $X$ is a $\bQ$-Gorenstein deformation of $Y$.
The minimal resolution of $\tilde{Y}$ is the toric surface which fibres over $\bP^1$ with two degenerate fibres, one of type
$(O)$ (where there are $m$ $(-2)$-curves as above) and one of type $(I)$ with self-intersection numbers
$$
-a_r,\ldots,-a_1,-(t+2),-b_1,\ldots,-b_{s+1},-1,-d_1,\ldots,-d_u,
$$ 
and two disjoint torus-invariant sections with self-intersection numbers $-a_{r+1}$ and $-b_{s+1}$ 
which intersect the first fibre in the two $(-1)$-curves and the second in the end components labelled $-a_r$ and $-d_u$ respectively.
Note that the sequence $d_1,\ldots,d_u$ is uniquely determined. Note also that, as above, the equality $m=a_{r+1}+b_{s+1}-2$
ensures that this does define a toric surface (it is obtained as a toric blowup of a $\bP^1$-bundle).
The surface $Y$ has singularities an $A_m$ singularity and the cyclic quotient singularities obtained by contracting the chains 
of smooth rational curves with self-intersection numbers $-a_{r+1},\ldots, -a_1$, $-(t+2),-b_1,\ldots,-b_{s+1}$ and 
$-d_1,\ldots,-d_u,-b_{s+1}$.
This last singularity is of type $T_{t+1}$ by Lem.~\ref{T}(2). 
Hence the singularities of $X$ are $\bQ$-Gorenstein deformations of the singularities of $Y$ --- 
the first two singularities are not deformed, and the $T_{t+1}$-singularity is deformed to an $A_t$ singularity. 
Moreover, this deformation does not change the Picard number. 
Let $X'$ be the general fibre of a $1$-parameter deformation of $X$ inducing this deformation of the singularities. 
We show that $X' \simeq X$.  

Let $\hat{Y} \rightarrow Y$ and $\hat{X}' \rightarrow X'$ be the minimal resolutions of the singularities we do not deform.
Thus $\hat{Y}$ has a single $T_{t+1}$ singularity given by contracting the chain of smooth rational curves with self-intersection numbers $-d_1,\ldots,-d_u,-b_{s+1}$ on $\tilde{Y}$. 
Let $C_1$ and $C_2$ be the images of the $(-1)$-curves on $\tilde{Y}$ incident to the ends of this chain.
Then $C_1$ and $C_2$ are smooth rational curves meeting in a node at the singular point.
We claim that $C=C_1+C_2$ deforms to a smooth $(-1)$-curve on $\hat{X}'$ (not passing through the singular point).
First, by Lem.~\ref{intersectionnumber} we have $C^2=-1$.
Second, we prove that $C$ deforms.
We work on the canonical covering stack $q \colon \hat{\cY} \rightarrow \hat{Y}$ of $\hat{Y}$, see Sec.~\ref{QG}. 
Note that the deformation of $\hat{Y}$ lifts to a deformation of $\hat{\cY}$ (because it is a $\bQ$-Gorenstein deformation).
Let $\cC \rightarrow C$ be the restriction of the covering $\hat{\cY} \rightarrow \hat{Y}$.
The closed substack $\cC \subset \hat{\cY}$ is a Cartier divisor.
Hence the obstruction to deforming $\cC \subset \hat{\cY}$ lies in $H^1(\cN_{\cC/\hat{\cY}})$, 
where $\cN_{\cC/\hat{\cY}}$ is the normal bundle $\cO_{\hat{\cY}}({\cC})|_{\cC}$.
We compute that this obstruction group is zero.
Consider the exact sequence
$$
0 \rightarrow \cN_{\cC/\hat{\cY}} \rightarrow \oplus \, \cN_{\cC/\hat{\cY}}|_{\cC_i} \rightarrow \cN_{\cC/\hat{\cY}} \otimes k(Q) \rightarrow 0.
$$
where $\cC_i \rightarrow C_i$ are the restrictions of $q$ and $Q \in \hat{\cY}$ is the point over the singular point $P \in \hat{Y}$.
Now push forward to the coarse moduli space $\hat{Y}$.
(Recall that if $\cX$ is a Deligne--Mumford stack and $q \colon \cX \rightarrow X$ is the map to its coarse moduli space, then locally over $X$ the map $q$ is of the form $[U/G] \rightarrow U/G$ where $U$ is a scheme and $G$ is a finite group acting on $U$. A sheaf $\cF$ over $[U/G]$ corresponds to a $G$-equivariant sheaf $\cF_U$ over $U$, and 
$q_*\cF=(\pi_*\cF_U)^G$ where $\pi \colon U \rightarrow U/G$ is the quotient map.)
Let $n$ be the index of the singularity $P \in Y$. Then $n>1$ and the $\bmu_n$ action on $\cN_{\cC/\hat{\cY}} \otimes k(Q)$ is non-trivial.
So $q_*(\cN_{\cC/\hat{\cY}} \otimes k(Q))=0$ and $q_*\cN_{\cC/\hat{\cY}}=\oplus \, q_*\cN_{\cC/\hat{\cY}}|_{\cC_i}$ by the exact sequence above.
The sheaf $q_*\cN_{\cC/\hat{\cY}}|_{\cC_i}$ is a line bundle on $C_i \simeq \bP^1$ of degree $\lfloor C \cdot C_i \rfloor$.
Let $\alpha \colon \tilde{Y} \rightarrow \hat{Y}$ denote the minimal resolution of $\hat{Y}$ and $C_i'$ the strict transform of $C_i$ for each $i$.
Then
$$C \cdot C_i = \alpha^*C \cdot C_i' > C_i'^2=-1.$$
Hence $H^1(q_*\cN_{\cC/\hat{\cY}}|_{\cC_i})=0$. We deduce that $H^1(\cN_{\cC/\hat{\cY}})=0$ as required.

We now compute locally that $C$ deforms to a smooth curve that does not pass through the singular point of $\hat{X}'$.
Locally at the singular point of $\hat{Y}$, the deformation of $\hat{Y}$ is of the form
$$(xy=(z^n-w)^d) \subset \frac{1}{n}(1,-1,a) \times \bC^1_w$$
where $d=t+1$.
The deformation of $C$ is given by an equation $(z+ w \cdot h = 0)$, where $h \in k[[x,y,w]]$ has $\bmu_n$-weight $a$.
So, eliminating $z$, the abstract deformation of $C$ is given by $(xy=u \cdot w^d) \subset \frac{1}{n}(1,-1) \times \bC^1_w$,
where $u$ is a unit. In particular the general fibre is smooth and misses the singular point of the ambient surface $\hat{X}'$.

We deduce that, on $\hat{X}'$, we have a cycle of smooth rational curves of self-intersections
$$-a_{r+1},\ldots,-a_1,-(t+2),-b_1,\ldots,-b_{s+1},-1,-2,\ldots,-2,-1$$
(where the chain of $(-2)$-curves has length $m$).
Indeed the chains $-a_{r+1}, \ldots, -b_{s+1}$ and $-2, \ldots, -2$ are the exceptional loci of the minimal resolutions of two of the singular points of $X'$, 
the first $(-1)$-curve is the deformation of $C$ described above, 
and the last $(-1)$-curve is the deformation of the $(-1)$-curve on $\hat{Y}$.
Moreover $\hat{X}'$ has a unique singular point of type $A_t$ which does not lie on this cycle.
Let $\tilde{X}' \rightarrow \hat{X}'$ be the minimal resolution. 
Observe that the chain $-1,-2,\ldots,-2,-1$ defines a ruling of $\tilde{X}'$. 
If $f$ is another degenerate fibre, then $f$ contains a unique $(-1)$-curve and its remaining components are exceptional over $X'$ (because $\rho(X')=1$). 
We deduce that there is exactly one additional degenerate fibre, which is the union of the chain $-a_r,\ldots,-b_s$, 
the chain $-2,\ldots,-2$ of length $t$ (the exceptional locus of the minimal resolution of the $A_t$ singularity) and a $(-1)$-curve. This determines the fibre uniquely.
We conclude that $X' \simeq X$.

A similar argument works when $\tilde{E}_1 \cdot \tilde{E}_2 = 1$.
In this case the ruling $p \colon \tilde{X} \rightarrow \bP^1$ has a unique degenerate fibre of type $(O)$
and the two sections $\tilde{E_1}$ and $\tilde{E_2}$ meet this fibre in the two $(-1)$-curves. 
Set $\tilde{E}_1^2=-a$ and $\tilde{E}_2^2=-b$ and let $m$ be the number of $(-2)$-curves in degenerate fibre.
Then $X$ has singularities $A_m$ and the cyclic quotient singularity whose minimal resolution has exceptional locus 
$\tilde{E_1}+\tilde{E_2}$.
(In particular, $(a,b)=(2,2)$, $(3,3)$, or $(2,5)$ because $X$ has $T$-singularities, but we give a uniform treatment of these
cases.)
We compute that $m=a+b+1$ by expressing $p$ as a blowup of a $\bP^1$-bundle.

We construct a toric surface $Y$ and prove that $X$ is a $\bQ$-Gorenstein deformation of $Y$.
The minimal resolution of $\tilde{Y}$ is the toric surface which fibres over $\bP^1$ with two degenerate fibres, one of type
$(O)$ (where there are $m$ $(-2)$-curves as above) and one of type $(I)$ with self-intersection numbers
$$
-b,-1,-2,\ldots,-2
$$ 
(where the chain of $(-2)$-curves has length $(b-1)$)
and two disjoint torus-invariant sections with self-intersection numbers $-a$ and $-(b+3)$ 
which intersect the first fibre in the two $(-1)$-curves and the second in the end components with self-intersection numbers 
$-b$ and $-2$ respectively.
Note that the equality $m=a+(b+3)-2$ ensures that this does define a toric surface.
The surface $Y$ has singularities an $A_m$ singularity and the cyclic quotient singularities obtained by contracting the chains 
of smooth rational curves with self-intersection numbers $-a,-b$ and 
$-2,\ldots,-2,-(b+3)$.
This last singularity is of type $T_{1}$ by Prop.~\ref{exclocusT_d}. 
Hence the singularities of $X$ are deformations of the singularities of $Y$ --- 
the first two singularities are not deformed, and the $T_1$-singularity is smoothed.
Moreover, this deformation does not change the Picard number. 
Let $X'$ be the general fibre of a $1$-parameter deformation of $X$ inducing this deformation of the singularities.  
Let $\hat{Y} \rightarrow Y$ and $\hat{X}' \rightarrow X'$ be the minimal resolutions of the singularities we do not deform.
Thus $\hat{Y}$ has a single $T_{1}$ singularity given by contracting the chain of smooth rational curves with self-intersection numbers $-2,\ldots,-2,-(b+3)$ on $\tilde{Y}$. 
Let $C_1$ and $C_2$ be the images of the $(-1)$-curves on $\tilde{Y}$ incident to the ends of this chain,
so $C_1$ and $C_2$ are smooth rational curves meeting in a node at the singular point.
Then, as above, $C=C_1+C_2$ deforms to a smooth $(-1)$-curve on $\hat{X}'$.
We deduce that, on $\hat{X}'$, we have a cycle of smooth rational curves of self-intersections
$$-a,-b,-1,-2,\ldots,-2,-1$$
(where the chain of $(-2)$-curves has length $m$).
Indeed, the chains $-a,-b$ and $-2, \ldots, -2$ are the exceptional loci of the minimal resolutions of the two singular points of 
$X'$, the first $(-1)$-curve is the deformation of $C$, and the last $(-1)$-curve is the deformation of the $(-1)$-curve on $\hat{Y}$.
Let $\tilde{X}' \rightarrow \hat{X}'$ be the minimal resolution. 
Observe that the chain $-1,-2,\ldots,-2,-1$ defines a ruling of $\tilde{X}'$. 
There are no other degenerate fibres of this ruling because $\rho(X')=1$.  
We deduce that $X' \simeq X$.

If we are in case $(2)$ of Thm.~\ref{fibration}, then the non Du Val singularities of $X$ are of type $\frac{1}{4}(1,1)$.
In particular, $2K_X$ is Cartier.
Similarly, if $X$ has a $D$ or $E$ singularity then $2K_X$ is Cartier by Thm.~\ref{DorE}. 
So in these cases we can refer to the classification of log del Pezzo surfaces of Picard rank $1$ and index $\le 2$ 
given by Alexeev and Nikulin \cite[Thms.~4.2,4.3]{AN}.
By Prop.~\ref{index2deformation} the only such surfaces which are not $\bQ$-Gorenstein deformations of toric surfaces
are those which have either a $D$ singularity, an $E$ singularity, or at least $4$ Du Val singularities.
These are the sporadic surfaces listed in Ex.~\ref{sporadic}. This completes the proof.
\end{proof}

\begin{lem}\label{T}
Let $[a_1,\ldots,a_r]$ and $[b_1,\ldots,b_s]$ be conjugate strings.
\begin{enumerate}
\item The conjugate of $[a_r,\ldots,a_1,t+2,b_1,\ldots,b_s]$ is a $T_{t+1}$-string.
\item Given $b_{s+1} \ge 2$, let $[d_1,\ldots,d_u]$ be the conjugate of $[a_r,\ldots,a_1,t+2,b_1,\ldots,b_s,b_{s+1}]$. 
Then $[d_1,\ldots,d_u,b_{s+1}]$ is a $T_{t+1}$-string.
\end{enumerate}
\end{lem}
\begin{proof}
Let an \emph{$S_t$-string} be a string $[a_r,\ldots,a_1,t+2,b_1,\ldots,b_s]$ as above. Then, by Lem.~\ref{conjugate},
we have
\begin{enumerate}
\item[$(a)$] $[2,t+2,2]$ is an $S_t$-string.
\item[$(b)$] If $[e_1,\ldots,e_v]$ is an $S_t$-string, then so are $[e_1+1,\cdots,e_v,2]$ and $[2,e_1,\ldots,$ $e_v+1]$.
\item[$(c)$] Every $S_t$-string is obtained from the example in 
$(a)$ by iterating the steps in $(b)$.
\end{enumerate}
Now $(1)$ follows from Prop.~\ref{exclocusT_d} and Lem.~\ref{conjugate}.
To deduce $(2)$, let $[e_1,\ldots,e_v]$ be the conjugate of $[a_r,\ldots,a_1,t+2,b_1,\ldots,b_s]$.
Then $$[d_1,\ldots,d_u,b_{s+1}]=[2,\ldots,2,e_1+1,e_2,\ldots,e_v,b_{s+1}]$$
(where there are $(b_{s+1}-2)$ 2's)  by Lem.~\ref{conjugate}. 
This string is of type $T_{t+1}$ by (1) and Prop.~\ref{exclocusT_d}.
\end{proof}

\begin{lem}\label{intersectionnumber}
Let $(P \in S,D)$ denote the local pair $(\frac{1}{dn^2}(1,dna-1),(uv=0))$.
Let $\pi \colon \tilde{S} \rightarrow S$ be the minimal resolution of $S$ and $D'$ the strict transform of 
$D$. Write $\pi^*D = D' + F$ where $F$ is $\pi$-exceptional. Then $F^2=-1$.
\end{lem}
\begin{proof}
We may assume $S$ is a projective toric surface, $P \in S$ is the unique singular point, and $D$ is the toric boundary.
Then $\tilde{S}$ is toric with boundary $\tilde{D}:=D'+\sum E_i$, where $E_1,\ldots,E_r$ are the exceptional divisors of $\pi$.
In particular $D \in |-K_S|$ and $\tilde{D} \in |-K_{\tilde{S}}|$.
Since $P \in S$ is a $T_d$-singularity, by Prop.~\ref{noether} we have
$$K_{\tilde{S}}^2+\rho(\tilde{S})=K_S^2+\rho(S)+(d-1).$$
So $\tilde{D}^2+r=D^2+(d-1)$. Now $\tilde{D}^2=D'^2+\sum E_i^2 +2(r+1)$, so 
$$F^2=D'^2-D^2=d-3r-3-\sum E_i^2.$$
Finally, $\sum E_i^2 = d-3r-2$ by the inductive description of resolutions of $T_d$-singularities (see Prop.~\ref{exclocusT_d}),
so $F^2=-1$ as claimed.
\end{proof} 

\begin{proof}[Proof of Thm.~{\ref{thm-ge}}]
Let $X$ denote the special fibre of $f \colon V \rightarrow T$.
Thus $X$ is a del Pezzo surface with quotient singularities which admits a $\bQ$-Gorenstein smoothing.
Since $H^1(\cO_X)=H^2(\cO_X)=0$ the restriction map $\Pic V \rightarrow \Pic X$
is an isomorphism. Hence $\rho(X)=\rho(V/T)=1$.

The pair $(V,X)$ is plt by inversion of adjunction \cite[Thm. 17.6]{FA}.
By Thm. \ref{1complement} there exists a 1-complement $D$ of $K_X$.
By \cite[Prop.~4.4.1]{P} $D$ lifts to a 1-complement $S$ of $K_V+X$. 
That is, there exists an effective divisor $S$ on $V$ such that $S|_X=D$, $K_V+X+S\sim 0$ (equivalently, $S\in |-K_V|$), and the pair $(V,X+S)$ is lc. 
It follows that the pair $(V,S)$ is also lc and has no log canonical centers contained in $X$.
A general fibre $(V_t,S_t)$ of $(V,S)/T$ is a smooth del Pezzo surface with anticanonical divisor.
So $S_t$ is smooth for $S$ general. We deduce that the pair $(V,S)$ is plt.
Thus $S$ is plt by adjunction, and Gorenstein, so has only Du Val singularities.
\end{proof}

Finally, we note that Cor.~\ref{deg9} follows from Thm.~\ref{thmintro}.
Indeed, if $X$ is a surface with quotient singularities which admits a smoothing to the plane,
then $\rho(X)=1$, $-K_X$ is ample, and the smoothing is $\bQ$-Gorenstein by \cite[\S 1]{M1}.

\section{Exceptional bundles associated to degenerations}

In this section we connect our results with the theory of exceptional vector bundles on smooth del Pezzo surfaces
\cite{Ru},\cite{Ru2},\cite{Ru3},\cite{KO},\cite{KN}. Full details will appear elsewhere.

Let $Y$ be a smooth projective surface. A vector bundle $F$ on $Y$ is \emph{exceptional} if $\End F = \bC$ and $\Ext^1(F,F)=\Ext^2(F,F)=0$.

\begin{thm}
Let $X$ be a projective surface with quotient singularities and $\cX/(0 \in T)$ a $\bQ$-Gorenstein smoothing.
Assume $H^1(\cO_X)=H^2(\cO_X)=0$ and $H_1(X^0,\bZ)=0$, where $X^0 \subset X$ is the smooth locus.
Let $P \in X$ be a singularity of type $\frac{1}{n^2}(1,na-1)$.
Then there exists a base change $T' \rightarrow T$ and a reflexive sheaf $\cE$ over $\cX'$ such that
\begin{enumerate}
\item $\cE|_{\cX'_t}$ is an exceptional vector bundle of rank $n$ on $\cX'_t$ for $0 \neq t \in T'$.
\item $E:=\cE|_X$ is a torsion-free sheaf on $X$ and there is an exact sequence
$$0 \rightarrow E \rightarrow \cO_X(D)^{\oplus n} \rightarrow \cT \rightarrow 0$$
where $D$ is a Weil divisor on $X$ and $\cT$ is a torsion sheaf supported at $P \in X$.
\end{enumerate}
\end{thm}

Fix $0 \neq t \in T'$ and write $Y=\cX'_t$ and $F=\cE|_Y$.
We say $F$ is a \emph{vanishing bundle} on $Y$ associated to $P \in X$.
The bundle $F$ is determined up to $F \mapsto F^{\vee}$ and $F \mapsto F \otimes L$ for $L \in \Pic Y$.

Note that there are no vanishing cycles for a $\bQ$-Gorenstein smoothing  of a $\frac{1}{n^2}(1,na-1)$ singularity by Lem.~\ref{milnornumber}.
We think of $F$ as analogous to a vanishing cycle.

A sequence $(F_0,\ldots,F_r)$ of exceptional bundles on a smooth surface $Y$ is an \emph{exceptional collection}
if $\Ext^i(F_j,F_k)=0$ for all $i$ and $j>k$.

\begin{thm}
With notation and assumptions as above, let $P_i \in X$ be a singularity of type $\frac{1}{n_i^2}(1,n_ia_i-1)$ for $i=0,\ldots,r$.
Let $\Gamma_1,\ldots,\Gamma_r$ be a chain of smooth rational curves connecting the $P_i$ such that, in orbifold coordinates $u_i,v_i$ at $P_i$,
we have $\Gamma_i=(v_i=0)$ and $\Gamma_{i+1}=(u_i=0)$. Then there exist vanishing bundles $F_i$ on $Y$ associated to $P_i \in X$ such that
$(F_0,\ldots,F_r)$ is an exceptional collection on $Y$.
\end{thm}

We say two exceptional collections are \emph{equivalent} if they are related by a sequence of operations of the following types.
\begin{enumerate}
\item $(F_0,\ldots,F_r) \mapsto (F_0 \otimes L,\ldots,F_r \otimes L)$, some $L \in \Pic Y$
\item $(F_0,\ldots,F_r) \mapsto (F_r^{\vee},\ldots,F_0^{\vee})$
\item $(F_0,\ldots,F_r) \mapsto (F_1,\ldots,F_r,F_0(-K_Y))$
\end{enumerate}

Combining Cor.~\ref{deg9} with the classification of exceptional collections on $\bP^2$ \cite{Ru}, we obtain

\begin{thm} There is a bijective correspondence between equivalence classes of exceptional collections on $\bP^2$ consisting of bundles of rank $>1$
and degenerations $X$ of $\bP^2$ with quotient singularities, given by the vanishing bundles.
\end{thm}

Suppose now that $(P \in X)$ is a singularity of type $\frac{1}{dn^2}(1,dna-1)$ and $(P \in \cX)/(0 \in T)$ is a $\bQ$-Gorenstein smoothing.
Then, after a base change $T' \rightarrow T$, there exists a simultaneous partial resolution $\pi \colon \tilde{\cX} \rightarrow \cX'/T'$ such that
$\pi_0 \colon \tilde{X} \rightarrow X$ has exceptional locus a chain of $(d-1)$ smooth rational curves connecting $d$ singularities of type $\frac{1}{n^2}(1,na-1)$, and $\pi_t$ is an isomorphism for $t \neq 0$. See \cite[Sec.~2]{BC}. So we can reduce to the case $d=1$ treated above.

Now, using Thm.~1.1, we can show that the three block complete exceptional collections on del Pezzo surfaces studied in \cite{KN}
arise from degenerations $X$ with $\rho(X)=1$.

\end{document}